\documentclass[final,3p,times]{elsarticle}

\usepackage{amssymb}
\usepackage[utf8]{inputenc}
\usepackage{amsmath}
\usepackage{amsfonts}
\usepackage{amsthm}
\usepackage{enumerate}
\usepackage{mathtools}
\usepackage{chngcntr} \usepackage{amsthm}

\theoremstyle{plain}
\newtheorem{theorem}{Theorem}[section]
\newtheorem{proposition}[theorem]{Proposition}
\newtheorem{lemma}[theorem]{Lemma}
\newtheorem{corollary}[theorem]{Corollary}

\theoremstyle{definition}

\newtheorem{definition}[theorem]{Definition}

\theoremstyle{remark}
\newtheorem{remark}[theorem]{Remark}

\counterwithin{equation}{section}

\providecommand{\abs}[1]{\lvert#1\rvert}
\providecommand{\norm}[1]{\lVert#1\rVert}

\journal{Journal de Mathématiques Pures et Appliquées}

\begin{document}

\begin{frontmatter}



\title{Leaves decompositions in Euclidean spaces}


\author{Krzysztof J. Ciosmak\corref{cor1}}
\cortext[cor1]{Funding: the financial support of St John’s College in Oxford, Clarendon Fund and EPSRC is gratefully acknowledged.  Part of this research was completed in Fall 2017 while the author was member of the Geometric Functional Analysis and Application program at MSRI, supported by the National Science Foundation under Grant No. 1440140. This research was also partly supported by the ERC Starting Grant 802689 CURVATURE.\\ Declaration of interest: none.}
\address{University of Oxford, Mathematical Institute,
Andrew Wiles Building, Radcliffe Observatory Quarter, Woodstock Rd, Oxford OX2 6GG, United Kingdom, \ead{ciosmak@maths.ox.ac.uk}\\
University of Oxford, St John’s College,
St Giles’, Oxford OX1 3JP,United Kingdom.
\ead{krzysztof.ciosmak@sjc.ox.ac.uk}}

\begin{abstract}

We partly extend the localisation technique from convex geometry to the multiple constraints setting.

For a given $1$-Lipschitz map $u\colon\mathbb{R}^n\to\mathbb{R}^m$, $m\leq n$, we define and prove the existence of a partition of $\mathbb{R}^n$, up to a set of Lebesgue measure zero, into maximal closed convex sets such that restriction of $u$ is an isometry on these sets. 

We consider a disintegration, with respect to this partition, of a log-concave measure. We prove that for almost every set of the partition of dimension $m$, the associated conditional measure is log-concave. This result is proven also in the context of the curvature-dimension condition for weighted Riemannian manifolds. This partially confirms a conjecture of Klartag.

\begin{flushleft} \bf{R\'esum\'e}
\end{flushleft}

Nous \'etendons en partie la technique de localisation de la g\'eom\'etrie convexe pour plusieurs contraintes.

\'Etant donn\'ee application $1$-lipschitzienne  $u\colon\mathbb{R}^n\to\mathbb{R}^m$, $m\leq n$, nous d\'efinissons et prouvons l'existence d'une partition de $\mathbb {R}^n$, en dehors d'un ensemble n\'egligeable, \`a ensembles convexes ferm\'es maximaux tels que la restriction de $u$ soit une isom\'etrie sur ces ensembles.

On consid\`ere une d\'esint\'egration, relativement \`a cette partition, d'une mesure log-concave. Nous montrons que pour presque chaque ensemble $m$-dimensionnelle de la partition, la mesure conditionnelle associ\'ee est log-concave. Ce r\'esultat est \'egalement prouv\'e dans le contexte de la condition de courbure-dimension pour les vari\'et\'es riemanniennes pond\'er\'ees. Cela confirme en partie une conjecture de Klartag.
\end{abstract}


\end{frontmatter}
\section{Introduction}

In this note we consider topics generalising the localisation technique and stemming from the optimal transport theory. Let us describe these connections.

\subsection{Optimal transport}

In 1781 Gaspard Monge (see \cite{Monge}) asked the following question: given two probability distributions $\mu,\nu$ on a metric space $(X,d)$, how to transfer one distribution onto the other in an optimal way. The criterion of optimality was to minimise the average transported distance. Since then the topic has been developed extensively and much of this development has been done recently. We refer the reader to the books of Villani (see \cite{Villani2}, \cite{Villani1}) and to the lecture notes of Ambrosio (see \cite{Ambrosio3}) for a thorough discussion, history and applications of the optimal transport problem. 

The modern mathematical treatment of the problem has been initiated in 1942 by Kantorovich \cite{KantorovichURSS}, \cite{Kantorovich}. He proposed to consider a relaxed problem of optimising
\begin{equation*}
\int_{X\times X}d(x,y)d\pi(x,y)
\end{equation*}
among all transference plans $\pi$ between $\mu$ and $\nu$, i.e., the set $\Pi(\mu,\nu)$ of Borel probability measures on $X\times X$ with respective marginal distributions equal to $\mu$ and to $\nu$.
The existence of an optimal transference plan is a straightforward consequence of the Prokhorov's theorem, provided that $X$ is separable.

The main question that has attracted a lot of attention is whether there exists an optimal transport plan, i.e., a Borel map $T\colon X\to X$ such that $T_{\#}\mu=\nu$ and the integral
\begin{equation*}
\int_Xd(x,T(x))d\mu(x)
\end{equation*}
is minimal.
If we knew that an optimal transference plan is concentrated on a graph of a Borel measurable function then we could infer the existence of an optimal transport plan. The first complete answer on Euclidean space, under regularity assumptions on the considered measures, was presented in a seminal paper \cite{Evans-Gangbo} of Evans and Gangbo. However, before that, Sudakov in \cite{Sudakov} presented a solution of the problem that contained a flaw. The flaw has been remedied by Ambrosio in \cite{Ambrosio3} and later by Trudinger and Wang in \cite{Trudinger} for the Euclidean distance and by Caffarelli, Feldman and McCann in \cite{Caffarelli} for distances induced by norms that satisfy certain smoothness and convexity assumptions. In \cite{Caravenna1} Caravenna has carried out the original strategy of Sudakov for general strictly convex norms and eventually Bianchini and Daneri in \cite{Bianchini1} accomplished the plan of a proof of Sudakov for general norms on finite-dimensional normed spaces. Let us note here also a paper \cite{Caravenna2}, which deals with a related problem in the context of faces of convex functions.

Let us describe briefly the strategy of Sudakov in the context of Euclidean spaces. We assume that the two Borel probability measures $\mu,\nu$ on $\mathbb{R}^n$ are absolutely continuous with respect to the Lebesgue measure. 

Let us recall the paramount Kantorovich--Rubinstein duality formula
\begin{equation*}
\sup\Big\{\int_{\mathbb{R}^n}u d(\mu-\nu)\mid u\text{ is }1\text{-Lipschitz}\Big\}=\inf\Big\{\int_{\mathbb{R}^n\times\mathbb{R}^n}\norm{x-y}d\pi(x,y)\mid \pi\in\Pi(\mu,\nu)\Big\}.
\end{equation*}
Let us take an optimal $u$ and an optimal $\pi$ in the two above optimisation problems. We may infer that 
\begin{equation*}
u(x)-u(y)=\norm{x-y}\text{ for }\pi\text{-almost every }(x,y)\in X\times X.
\end{equation*}
Consider the maximal sets on which $u$ is an isometry, called the \emph{transport rays}. We see that all transport has to occur on these sets. Careful analysis of the Lipschitz function $u$ shows that the transport rays form a foliation of the underlying space $\mathbb{R}^n$, up to Lebesgue measure zero. It turns out that the direction of the transport rays is itself locally Lipschitz. This allows us to use of the area formula, which yields that the conditional measures of the disintegration of the Lebesgue measure with respect to the aforementioned foliation are absolutely continuous with respect to the one-dimensional Hausdorff measures on the transport rays. This is exactly the place where Sudakov's proof in \cite{Sudakov} contained a defect. He claimed that any foliation into segments is such that the conditional measures are absolutely continuous with respect to the one-dimensional Hausdorff measure. It was later shown by Ambrosio, Kirchheim and Preiss (see \cite{Ambrosio2}) that there exists a foliation consisting of segments and an atomic distributions on each segment such that the averaged measure is absolutely continuous with respect to the Lebesgue measure, refuting the claim of Sudakov.

Knowing that the conditional measures are absolutely continuous with respect to the one-dimensional Hausdorff measure, we may apply the well understood one-dimensional theory, where an optimal transport plan is known to exist and may be given by a certain formula, provided that at least one of the measures is non-atomic.
Then the optimal transport plan on the whole space is defined separately on each transport ray. 

The ideas of Sudakov have been applied also to other settings than normed spaces. The strategy has been carried out also in the context of Riemannian manifolds by Feldman and McCann in \cite{Feldman}.

\subsection{Localisation technique}

In \cite{Klartag} Klartag has observed that the above described strategy of Sudakov, in its instances in works of Caffarelli, Feldman and McCann \cite{Caffarelli} and of Feldman and McCann \cite{Feldman}, may be applied to adapt the \emph{localisation technique} from convex geometry to the setting of Riemannian manifolds.
The technique allows to reduce certain high dimensional problems with one linear constraint to a collection of one-dimensional problems with an analogous constraint. 
Let us include a brief description of the technique based on \cite{Klartag}. 

It first appeared in works of Payne and Weinberger \cite{Payne} and was developed in the context of convex geometry by Gromov and Milman \cite{Gromov}, Lov\'asz and Simonovits \cite{Lovasz1} and by Kannan, Lov\'asz and Simonovits \cite{Lovasz2}. Later, Klartag \cite{Klartag} adapted the technique to the setting of weighted Riemannian manifolds satisfying the curvature-dimension condition in the sense of Bakry and \'Emery  \cite{Bakry1}, \cite{Bakry}. Subsequently, Ohta \cite{Ohta} generalised these results to Finsler manifolds and Cavalletti and Mondino \cite{Cavalletti3}, \cite{Cavalletti2} generalised them to metric measure spaces satisfying the synthetic curvature-dimension condition. The latter was introduced in the foundational papers by Sturm \cite{Sturm1}, \cite{Sturm2} and by Lott and Villani \cite{Villani3} and allowed for development of a far-reaching, vast theory of metric measure spaces. The curvature-dimension condition may be thought of as lower bound on the curvature and an upper bound on the dimension of the considered space. We refer also to Ambrosio \cite{Ambrosio4} for a recent account on the spaces satisfying the curvature-dimension condition. Let us note that the curvature-dimension condition is also related to Bochner's inequality; see \cite{Sturm3}.

The technique developed by Payne and Weinberger has no clear analogue for an abstract Riemannian manifold. This is the point where the optimal transport plays its r\^ole in localisation. Let us cite below a theorem from \cite{Klartag}, presented there in a general setting of Riemannian manifolds and measures that satisfy the curvature-dimension condition; see Definition \ref{defin:curv}.

\begin{theorem}\label{thm:localisation}
Let $n\geq 2$, $\kappa\in\mathbb{R}$ and $N\in (-\infty,1)\cup[n,\infty]$. Assume that $(\mathcal{M},d,\mu)$ is a geodesically convex, $n$-dimensional weighted Riemannian manifold satisfying the curvature-dimension condition $CD(\kappa,N)$ w. Let $g\colon \mathcal{M}\to\mathbb{R}$ be a $\mu$-integrable function such that
\begin{equation*}
\int_{\mathcal{M}}gd\mu=0\text{ and }\int_{\mathcal{M}}\abs{g(x)}d(x,x_0)d\mu(x)<\infty\text{ for some }x_0\in\mathcal{M}.
\end{equation*}
Then there exists a partition $\Omega$ of $\mathcal{M}$ into pairwise disjoint sets, a measure $\nu$ on $\Omega$ and a family $(\mu_I)_{I\in\Omega}$ of measures on $\mathbb{R}^n$ such that:
\begin{enumerate}[i)]
\item\label{i:mixture} for any Lebesgue measurable set $A\subset\mathcal{M}$ the map $I\mapsto \mu_I(A)$ is well-defined $\nu$-almost everywhere, is $\nu$-measurable and
\begin{equation*}
\mu(A)=\int_{\Omega}\mu_I(A)d\nu(I),
\end{equation*}
\item\label{i:geodesic} for $\nu$-almost every $I\in\Omega$ the set $I\subset\mathcal{M}$ is a minimising geodesic and $\mu_I$ is supported on $I$ and is a $CD(\kappa,N)$-needle or else it is a singleton,
\item\label{i:balance} for $\nu$-almost every $I\in\Omega$ we have $\int_Igd\mu_I=0$.
\end{enumerate}
\end{theorem}

Let us remark that the above theorem has been known before in the context of Euclidean spaces. Let us note that the proof presented in \cite{Klartag} differs much from the previously known proofs of Gromov \cite{Gromov} or of Lov\'asz and Simonovits \cite{Lovasz1}, which employed the Borsuk--Ulam theorem.

The purposes of this article is to continue along this line of research and investigate multi-dimensional analogue of the localisation technique, as proposed in \cite[Chapter 6]{Klartag}, thus further generalising ideas of Sudakov.

To this end, we  shall consider finite-dimensional linear spaces equipped with Euclidean norm and $1$-Lipschitz map $u\colon\mathbb{R}^n\to\mathbb{R}^m$. In Section \ref{sec:partition} we define a partition  of $\mathbb{R}^n$, up to a set of Lebesgue measure zero, associated to such a map and prove its basic properties; see Lemma \ref{lem:diff}, Lemma \ref{lem:boundary}, Corollary \ref{col:unique}. The elements of the partition are maximal sets $\mathcal{S}$ such that the restriction of $u$ to $\mathcal{S}$ is an isometry, i.e. preserves the Euclidean distance. Each such set we shall call a \emph{leaf} of $u$. We prove that each leaf of $u$ is closed and convex (see Corollary \ref{col:convex}) hence it has a well-defined dimension. This is a multi-dimensional generalisation of \ref{i:geodesic}) of Theorem \ref{thm:localisation} and of the ideas concerning transport rays from optimal transport. Thanks to Lemma \ref{lem:important} and Corollary \ref{col:strength} we will provide a significantly simpler proof of Lemma \ref{lem:efge} than the proof of the analogous result in \cite{Caffarelli}.

In Theorem \ref{thm:dis} we show that we may decompose the Lebesgue measure on $\mathbb{R}^n$ into a mixture of measures, each supported on a leaf of $u$. In particular, the same is true for any measure $\mu$ such that $(\mathbb{R}^n,\norm{\cdot},\mu)$ satisfies the $CD(\kappa,N)$ condition, as any such measure is absolutely continuous with respect to the Lebesgue measure. It is a step towards a conjecture of Klartag \cite[Chapter 6]{Klartag} and a generalisation of \ref{i:mixture}) of Theorem \ref{thm:localisation}.

Suppose now that $m\leq n$ and that $(\mathbb{R}^n,\norm{\cdot},\mu)$ is a weighted Riemannian manifold, satisfying the curvature-dimension condition $CD(\kappa,N)$ for some $\kappa\in\mathbb{R}$ and $N\in (-\infty,1)\cup [n,\infty]$; see Definition \ref{defin:curv}. Here $\norm{\cdot}$ denotes the Euclidean metric on $\mathbb{R}^n$ and $\mu$ is a Borel finite measure on $\mathbb{R}^n$. 
A partial affirmative answer to the conjecture of Klartag is provided by Theorem \ref{thm:discd}, where we prove that, for the leaves $\mathcal{S}$ of $u$ of dimension $m$, the conditional measures $\mu_{\mathcal{S}}$ are supported on the relative interiors $\mathrm{int}\mathcal{S}$ and are such that $(\mathrm{int}\mathcal{S},\norm{\cdot},\mu_{\mathcal{S}})$ satisfies $CD(\kappa,N)$. This further developes the generalisation of \ref{i:geodesic}) of Theorem \ref{thm:localisation}.

Note that in \cite[Chapter 6]{Klartag}, it is conjectured that also the above theorem holds true also for leaves of $u$ of arbitrary dimension. 

Let us note that Theorem \ref{thm:discd}, described above, proves in particular that if we disintegrate the Lebesgue measure with respect to the partition obtained from a $1$-Lipschitz map, then $(\mathrm{int}\mathcal{S},\norm{\cdot},\mu_{\mathcal{S}})$ will satisfy the curvature-dimension condition $CD(0,n)$ for leaves $\mathcal{S}$ of dimension $m$. This complements the results of \cite{Ambrosio3}, \cite{Sudakov}; see also \cite{Caravenna1}, \cite{Caravenna2} and \cite{Bianchini1}. Note that our result tells in particular that the conditional measures are equivalent to the $m$-dimensional Hausdorff measure, which provides a strengthening of the previously known results. Note also that the condition $CD(0,\infty)$ is equivalent to log-concavity of the considered measure.

The possible applications of the results of the article are in the localisation or dimensional reduction arguments, where the disintegration is an effective tool. A similar result to ours in case $m=1$ has been used to derive new proofs and generalisations of isoperimetric inequality, Poincar\'e's inequality and others to the setting of metric measure spaces satisfying curvature bounds. We refer the reader to \cite{Cavalletti3}, \cite{Cavalletti2}, \cite{Klartag}, \cite{Ohta}.

The proof relies on the area formula and Fubini's theorem and is based on a work of Caffarelli, Feldman and McCann \cite{Caffarelli} and of Klartag \cite{Klartag}. See also \cite{Ambrosio3} and \cite{Feldman} for similar approach to the Monge--Kantorovich problem.

Another tool that we use is the Wijsman topology \cite{Wijsman} on the convex and closed subsets $CC(\mathbb{R}^n)$ of $\mathbb{R}^n$ which makes it a Polish space, so we may apply disintegration theorem. This use is inspired by a paper of Ob\l\'oj and Siorpares \cite{Siorpaes}.

Let us mention here a development \cite{Ciosmak1}, where a generalisation of optimal transport to vector measures is studied. In there, it is shown that the mass-balance condition, of vital importance for the classical optimal transport problem, does not hold for absolutely continuous vector measures, thus resolving another conjecture of Klartag \cite[Chapter 6]{Klartag} in the negative. This is to say, the proposed generalisation of \ref{i:balance}) of Theorem \ref{thm:localisation} does not hold true in the multiple-dimensional setting. Let us note that the outline of a proof from \cite{Klartag} of the conjecture has a gap, as follows by the results of \cite{Ciosmak}.

\subsection{Waist inequalities}

The waist inequality, proved by Gromov \cite{Gromov3}, \cite{Gromov2}, states that if $f\colon \mathbb{R}^n\to\mathbb{R}^m$, $m\leq n$ is a continuous function, then there exists $t\in\mathbb{R}^m$ such that the fibre $L=f^{-1}(t)$ satisfies
\begin{equation*}
\gamma_n(L+rB_n)\geq \gamma_m(rB_m)\text{ for all }r>0.
\end{equation*}
Here $\gamma_n$ and $\gamma_m$ are the $n$ and $m$ dimensional standard Gaussian measures, $B_n$ and $B_m$ are the unit balls in $\mathbb{R}^n$ and $\mathbb{R}^m$ respectively. This inequality may be seen as a generalisation of the Gaussian isoperimetric inequality. Gromov \cite{Gromov2} has provided a proof of this inequality with use of the localisation method \cite{Gromov3}, \cite{Lovasz2}, \cite{Lovasz1}, \cite{Payne} combined with a Borsuk--Ulam type theorem. Later, Klartag \cite{Klartag2} has proved the theorem for the unit cube also with use of localisation methods, confirming a conjecture of Guth \cite{Guth}. 

One of the future possible applications of the research initiated in this article is to prove a general version of the inequality for spaces satisfying the curvature-dimension condition. This would imply the version for convex bodies and, in turn, would help answering the Bourgain’s hyperplane conjecture \cite{Bourgain} and the isoperimetric conjecture of Kannan, Lov\'{a}sz and Simonovits \cite{Lovasz2}.

\subsection{Multi-bubble conjectures}

The Gaussian multi-bubble conjecture is a generalisation of isoperimetric inequality that states that among all decompositions of $\mathbb{R}^n$ into $2\leq k\leq n+1$ sets of prescribed Gaussian measure the minimal Gaussian-weighted perimeter is uniquely attained by the Voronoi cells of $k$ equidistant points. The conjecture has been recently confirmed by E. Milman and Neeman \cite{Milman2} (c.f. \cite{Milman1}). Another possible application of the multi-dimensional localisation is a generalisation of this inequality for spaces satisfying the curvature-dimension condition.

\subsection{Outiline of the article}

In Section \ref{sec:partition} we provide a definition of the partition associated to any $1$-Lipschitz map $u\colon\mathbb{R}^n\to\mathbb{R}^m$. We prove that certain components of $u$ are differentiable on certain leaves, see Lemma \ref{lem:diff}. Moreover we investigate the regularity of the derivative on the leaves and provide a strengthening of $1$-Lipschitz property of $u$; see Lemma \ref{lem:important} and Remark \ref{rmk:strength}.
See also Lemma \ref{lem:boundary}, Corollary \ref{col:unique} for results concerning disjointness of the elements of the partition.

In Section \ref{sec:varia} we define a Lipschitz change of variables on so-called clusters of leaves, that will allow us to use the area formula and then Fubini's theorem to prove the regularity properties of the conditional measures; see Lemma \ref{lem:efge}. 

In Section \ref{sec:measur} we prove measurability properties of the partition; see Corollary \ref{col:borel}. We also prove that the union of boundaries of leaves of maximal dimension is a Borel set of the Lebsegue measure zero; see Lemma \ref{lem:measurable}. The material included here concerning leaves of non-maximal dimension is not employed in further investigations.

In Section \ref{sec:disin} we provide a proof of Theorem \ref{thm:dis}, that the partition induces a disintegration of the Lebesgue measure.

In Section \ref{sec:curv} we prove that the weighted Riemannian manifolds $(\mathrm{int}\mathcal{S},d,\mu_{\mathcal{S}})$ satisfy the curvature-dimension condition, provided that $(\mathbb{R}^n, d, \mu)$ did; see Theorem \ref{thm:discd}. This partially resolves in the affirmative a conjecture of Klartag \cite[Chapter 6]{Klartag}.

\section{Partition and its regularity}\label{sec:partition}

If $A\subset \mathbb{R}^n$ let us denote by $\mathrm{Conv}A$ the \emph{convex hull} of $A$, i.e. 
\begin{equation*}
\mathrm{Conv}A=\Big\{\sum_{i=1}^k \lambda_i x_i\mid k\in\mathbb{N},\lambda_1,\dotsc,\lambda_k\geq 0, \sum_{i=1}^k\lambda_i=1,x_1,\dotsc,x_k\in A\Big\}.
\end{equation*}
We define the \emph{affine hull} $\mathrm{Aff}A$ of a set $A\subset\mathbb{R}^n$ to be 
\begin{equation*}
\mathrm{Aff}A=\Big\{\sum_{i=1}^k \lambda_i x_i\mid k\in\mathbb{N},\lambda_1,\dotsc,\lambda_k\in\mathbb{R}, \sum_{i=1}^k\lambda_i=1,x_1,\dotsc,x_k\in A\Big\}.
\end{equation*}

\begin{lemma}\label{lem:coordinates2}
Let $z_1,\dotsc,z_k\in\mathbb{R}^n$. 
Let $x,y\in \mathbb{R}^n$. Suppose that 
\begin{equation*}
\norm{x-z_i}\leq\norm{y-z_i},
\end{equation*}
for $i=1,\dotsc, k$. Then for all $z\in\mathrm{Conv}\{z_1,\dotsc.z_k\}$ there is
\begin{equation*}
\norm{x-z}\leq\norm{y-z}.
\end{equation*}
In particular, if $y\in\mathrm{Conv}\{z_1,\dotsc.z_k\}$, then $x=y$.
\end{lemma}
\begin{proof}
Let 
\begin{equation*}
Z=\mathrm{Conv}\{z_1,\dotsc,z_k\}.
\end{equation*}
We have
\begin{equation*}
\norm{x}^2+\norm{z_i}^2-2\langle x,z_i\rangle\leq \norm{y}^2+\norm{z_i}^2-2\langle y,z_i\rangle
\end{equation*}
for all $i=1,\dotsc,k$. Hence, for these $i$'s, we have 
\begin{equation*}
\norm{x}^2-2\langle x,z_i\rangle\leq\norm{y}^2-2\langle y,z_i\rangle
\end{equation*}
Thus, adding up these inequalities multiplied by non-negative coefficients that sum up to one, we get
\begin{equation*}
\norm{x}^2-2\langle x,z\rangle\leq\norm{y}^2-2\langle y,z\rangle
\end{equation*}
for all $z\in Z$. Hence also
\begin{equation*}
\norm{x-z}\leq\norm{y-z}.
\end{equation*}
Putting $z=y$ yields $\norm{x-y}=0$.
\end{proof}

Let $A\subset\mathbb{R}^n$. We shall say that a map $v\colon A\to\mathbb{R}^m$ is an isometry provided that for all $x,y\in A$ there is $\norm{v(x)-v(y)}=\norm{x-y}$.

\begin{definition}
Let $u\colon \mathbb{R}^n\to\mathbb{R}^m$ be a $1$-Lipschitz map. A set $\mathcal{S}\subset \mathbb{R}^n$ is called a \emph{leaf} of $u$ if $u|_{\mathcal{S}}$ is an isometry and for any $y\notin \mathcal{S}$ there exists $x\in\mathcal{S}$ such that $\norm{u(y)-u(x)}<\norm{y-x}$.
\end{definition}

In other words, $\mathcal{S}$ is a leaf if it is a maximal set, with respect to the order induced by inclusion, such that $u|_{\mathcal{S}}$ is an isometry.

\begin{definition}
If $C\subset\mathbb{R}^n$ is a convex set, then we shall call the \emph{tangent space} of $C$ the linear space $\mathrm{Aff}(C)-\mathrm{Aff}(C)$.
We shall call the \emph{relative interior} of $C$ the relative interior with respect to the topology of $\mathrm{Aff}(C)$.
\end{definition}

\begin{lemma}\label{lem:unique}
Let $\mathcal{S}\subset \mathbb{R}^n$ be an arbitrary subset. Let $u\colon \mathcal{S}\to \mathbb{R}^m$ be an isometry. Then there exists a unique $1$-Lipschitz function $\tilde{u}\colon \mathrm{Conv}(\mathcal{S})\to \mathbb{R}^m$ such that $\tilde{u}|_{\mathcal{S}}=u$. Moreover $\tilde{u}$ is an isometry.
\end{lemma}
\begin{proof}
Observe that, by the polarisation formula, $u$ preserves the scalar product, that is for all points $p,q,r,s\in \mathcal{S}$ there is
\begin{equation}\label{eqn:polarisation}
\begin{aligned}
&\langle u(p)-u(q), u(r)-u(s)\rangle =\\
&=\frac 12 \big( \norm{u(p)-u(s)}^2+\norm{u(q)-u(r)}^2-\norm{u(p)-u(r)}^2-\norm{u(q)-u(s)}^2\big)=\\
&=\frac 12 \big( \norm{p-s}^2+\norm{q-r}^2-\norm{p-r}^2-\norm{q-s}^2\big)=\langle p-q,r-s\rangle.
\end{aligned}
\end{equation}
Suppose that $y_1,\dotsc,y_k, z_1,\dotsc, z_l\in \mathcal{S}$ and that $s_1,\dotsc,s_k$, $t_1,\dotsc,t_l$ are non-negative real numbers such that
\begin{equation*}
\sum_{i=1}^ks_i=\sum_{j=1}^l t_j=1.
\end{equation*}
Then, by (\ref{eqn:polarisation}),
\begin{equation}\label{eqn:affinity}
\begin{aligned}
&\Big\lVert\sum_{i=1}^k s_iu(y_i)-\sum_{j=1}^lt_j u(z_j)\Big\rVert^2=\Big\lVert\sum_{i=1}^k\sum_{j=1}^l s_it_j(u(y_i)-u(z_j))\Big\rVert^2=\\
&=\sum_{i,i'=1}^k\sum_{j,j'=1}^l s_is_{i'}t_jt_{j'} \langle u(y_i)-u(z_j),u(y_{i'})-u(z_{j'})\rangle=\\
&=\sum_{i,i'=1}^k\sum_{j,j'=1}^l s_is_{i'}t_jt_{j'} \langle y_i-z_j,y_{i'}-z_{j'}\rangle=\Big\lVert\sum_{i=1}^k s_iy_i-\sum_{j=1}^lt_jz_j\Big\rVert^2.
\end{aligned}
\end{equation}
We may now affinely extend $u$ to $\mathrm{Conv}(\mathcal{S})$. That is, if $x_1,\dotsc,x_k\in \mathcal{S}$ and $s_1,\dotsc,s_k$ are any non-negative real numbers that sum up to one, we set
\begin{equation*}
\tilde{u}\bigg(\sum_{i=1}^k s_ix_i\bigg)=\sum_{i=1}^k s_iu(x_i).
\end{equation*}
Now, (\ref{eqn:affinity}) shows that $\tilde{u}$ is a well-defined affine map on $\mathrm{Conv}(\mathcal{S})$ and that it is an isometry.

Suppose now that we have another $1$-Lipschitz extension $v\colon \mathrm{Conv}(\mathcal{S})\to\mathbb{R}^m$. To prove that $v=\tilde{u}$ it is enough to show that $v$ is affine. Choose non-negative real numbers $s_1,\dotsc,s_k$ summing up to one and any points $x_1,\dotsc,x_k\in \mathcal{S}$.
Then, by $1$-Lipschitzness and by the fact that $v$ is isometric on $\mathcal{S}$, we get, as in (\ref{eqn:affinity}),
\begin{equation*}
\Big\lVert v\Big(\sum_{i=1}^ks_ix_i\Big)-v(x_j)\Big\rVert\leq \Big\lVert \sum_{i=1}^ks_ix_i-x_j\Big\rVert =\Big\lVert\sum_{i=1}^ks_iv(x_i)-v(x_j)\Big\rVert.
\end{equation*}
By Lemma \ref{lem:coordinates2} we see that 
\begin{equation*}
v\bigg(\sum_{i=1}^ks_ix_i\bigg)=\sum_{i=1}^ks_iv(x_i).
\end{equation*}
It follows that $v$ is affine on $\mathrm{Conv}(\mathcal{S})$.
\end{proof}

\begin{corollary}\label{col:convex}
Any leaf $\mathcal{S}$ of $u$ is a closed convex set and $u|_{\mathcal{S}}$ is an affine isometry.
\end{corollary}

Let $\mathcal{S}$ be a leaf of $u$. Let $P$ denote the orthogonal projection of $\mathbb{R}^n$ onto the tangent space $V$ of $\mathcal{S}$. Let 
\begin{equation*}
T\colon V\to\mathbb{R}^m
\end{equation*}
be a linear isometry such that 
\begin{equation*}
u(x)-u(y)=T(x-y)
\end{equation*}
for any $x,y\in\mathcal{S}$. Let $Q$ denote the orthogonal projection of $\mathbb{R}^m$ onto $T(V)$.

Below by $\mathrm{int}\mathcal{S}$, $\mathrm{cl}\mathcal{S}$, $\partial\mathcal{S}$ we understand the relative interior, the relative closure and the relative boundary of $\mathcal{S}$ respectively.

\begin{lemma}\label{lem:important}
Let $u\colon \mathbb{R}^n\to\mathbb{R}^m$ be a $1$-Lipschitz map. Let $\mathcal{S}_1,\mathcal{S}_2$ be two leaves of $u$. Let $V_1,V_2$ be their respective tangent spaces and let $P_1,P_2$ be orthogonal projections onto $V_1,V_2$ respectively. Let $T_1,T_2$ be isometric maps such that 
\begin{equation*}
u(x)-u(y)=T_i(x-y)\text{ for all }x,y\in\mathcal{S}_i, i=1,2.
\end{equation*}
Let $x_i\in\mathcal{S}_i$ and $\sigma_i=\mathrm{dist}(x_i,\partial \mathcal{S}_i)$ for $i=1,2$. Then
\begin{equation*}
2\sigma_1\sigma_2\norm{P_1P_2-P_1T_1^*T_2P_2}\leq \norm{x_1-x_2}^2-\norm{u(x_1)-u(x_2)}^2,
\end{equation*}
and for $i=1,2$
\begin{equation*}
2\sigma_i\norm{P_iT_i^*(u(x_1)-u(x_2))-P_i(x_1-x_2)}\leq\norm{x_1-x_2}^2-\norm{u(x_1)-u(x_2)}^2.
\end{equation*}
\end{lemma}
\begin{proof}
Let $y_i\in\mathcal{S}_i$ for $i=1,2$. Let $v_i=y_i-x_i$ for $i=1,2$. Then we may write
\begin{equation*} 
u(y_1)-u(y_2)=u(x_1)-u(x_2)+T_1v_1-T_2v_2.
\end{equation*}
Hence $\norm{u(y_1)-u(y_2)}^2$ is equal to 
\begin{equation*}
\norm{u(x_1)-u(x_2)}^2+\norm{v_1}^2+\norm{v_2}^2+2\langle u(x_1)-u(x_2),T_1v_1-T_2v_2\rangle-2\langle T_1v_1,T_2v_2\rangle.
\end{equation*}
We also have
\begin{equation*}
y_1-y_2=x_1-x_2+v_1-v_2,
\end{equation*}
yielding
\begin{equation*}
\norm{y_1-y_2}^2=\norm{x_1-x_2}^2+\norm{v_1}^2+\norm{v_2}^2+2\langle x_1-x_2,v_1-v_2\rangle-2\langle v_1,v_2\rangle.
\end{equation*}
As $u$ is $1$-Lipschitz, $\norm{u(y_1)-u(y_2)}\leq\norm{y_1-y_2}$. By the two identities above we get therefore that 
\begin{equation*}
2\langle v_1,v_2\rangle-2\langle T_1v_1,T_2v_2\rangle+2\langle u(x_1)-u(x_2),T_1v_1-T_2v_2\rangle-2\langle x_1-x_2,v_1-v_2\rangle
\end{equation*}
is bounded above by 
\begin{equation*}
\norm{x_1-x_2}^2-\norm{u(x_1)-u(x_2)}^2.
\end{equation*}
Suppose that $\sigma_1, \sigma_2$ are both positive. As $y_1,y_2$ were arbitrary points of $\mathcal{S}_1,\mathcal{S}_2$ respectively, the above inequality holds true for any $v_1\in V_1$ and any $v_2\in V_2$ of norm at most $\sigma_1$ and $\sigma_2$ respectively. If we add two such inequalities with $v_1,v_2$ replaced by $-v_1,-v_2$ then we get that 
\begin{equation}\label{eqn:inequalityimportant}
2\langle v_1,v_2\rangle-2\langle T_1v_1,T_2v_2\rangle\leq \norm{x_1-x_2}^2-\norm{u(x_1)-u(x_2)}^2.
\end{equation}
Equivalently, for any $w_1,w_2\in\mathbb{R}^n$ of norm at most one, we have
\begin{equation*}
2\sigma_1\sigma_2\big\langle w_1,(P_1P_2-P_1T_1^*T_2P_2)w_2\big\rangle\leq \norm{x_1-x_2}^2-\norm{u(x_1)-u(x_2)}^2.
\end{equation*}
Taking supremum over all $w_1,w_2\in\mathbb{R}^n$ of norm at most one yields the first desired inequality.
For the next inequalities, we assume that $\sigma_2>0$ and we take $v_1=0$ to get that
\begin{equation*}
-2\langle u(x_1)-u(x_2),T_2v_2\rangle+2\langle x_1-x_2,v_2\rangle\leq \norm{x_1-x_2}^2-\norm{u(x_1)-u(x_2)}^2.
\end{equation*}
Analogously for $v_2=0$ and $\sigma_1>0$
\begin{equation*}
2\langle u(x_1)-u(x_2),T_1v_1\rangle-2\langle x_1-x_2,v_1\rangle\leq \norm{x_1-x_2}^2-\norm{u(x_1)-u(x_2)}^2.
\end{equation*}
Hence for any $w_1,w_2\in\mathbb{R}^n$ of norm at most one there is
\begin{equation*}
\sigma_2\Big\langle\big( P_2T_2^*\big(u(x_1)-u(x_2)\big)-P_2(x_1-x_2)\big),w_2\Big\rangle\leq \norm{x_1-x_2}^2-\norm{u(x_1)-u(x_2)}^2
\end{equation*}
and
\begin{equation*}
\sigma_1\Big\langle\big( P_1T_1^*\big(u(x_1)-u(x_2)\big)-P_1(x_1-x_2)\big),w_1\Big\rangle\leq \norm{x_1-x_2}^2-\norm{u(x_1)-u(x_2)}^2.
\end{equation*}
Taking suprema over $w_1,w_2$ in the unit ball of $\mathbb{R}^n$ yields the desired results.
\end{proof}

\begin{remark}\label{rmk:strength}
Lemma \ref{lem:important} tells us that if $x_1,x_2$ belong to relative interiors of leaves $\mathcal{S}_1,\mathcal{S}_2$ respectively, then the $1$-Lipschitzness of map $u\colon\mathbb{R}^n\to\mathbb{R}^m$ is strengthened to the condition that
\begin{equation*}
\norm{u(x_1)-u(x_2)}^2+2\sigma_1\sigma_2\norm{P_1P_2-P_1T_1^*T_2P_2}\leq \norm{x_1-x_2}^2
\end{equation*} 
for all $x_1\in\mathcal{S}_1$ and all $x_2\in\mathcal{S}_2$.
\end{remark}

\begin{lemma}\label{lem:diff}
Let $\mathcal{S}$ be a leaf of a $1$-Lipschitz map $u\colon \mathbb{R}^n\to\mathbb{R}^m$. Then $Qu$ is differentiable in the relative interior of $\mathcal{S}$. Moreover, if $z_0$ belongs to the relative interior of $\mathcal{S}$, then
\begin{equation*}
DQu(z_0)=TP.
\end{equation*}
If $u$ is differentiable in $z_0$ for some $z_0\in\mathcal{S}$, then
\begin{equation*}
QDu(z_0)=TP.
\end{equation*}
\end{lemma}
\begin{proof}
Observe that $Q=TT^*$. Hence, by Lemma \ref{lem:important}, we see that 
\begin{equation*}
2\sigma \norm{Q (u(z_1)-u(z_0))-TP(z_1-z_0)}\leq \norm{z_1-z_0}^2-\norm{u(z_1)-u(z_0)}^2.
\end{equation*}
for all $z_0\in\mathcal{S}$ and $z_1\in\mathbb{R}^n$. Here $\sigma=\mathrm{dist}(z_0,\partial\mathcal{S})$. Hence if $\sigma>0$ we obtain that
\begin{equation*}
\limsup_{z_1\to z_0} \frac{\norm{Q(u(z_1)-u(z_0))-TP(z_1-z_0)}}{\norm{z_1-z_0}}\leq \limsup_{z_1\to z_0} \frac{\norm{z_1-z_0}}{2\sigma}=0.
\end{equation*}
This yields the asserted differentiability. 

Now, suppose that $u$ is differentiable at $z_0\in\mathcal{S}$.  Arguing as in the proof of Lemma \ref{lem:important} we see that for all $z_2\in \mathcal{S}$ and $z_1\in\mathbb{R}^n$ we have
\begin{equation*}
2\big\langle T^*(u(z_0)-u(z_1))-(z_0-z_1),z_2-z_0\big\rangle\leq \norm{z_1-z_0}^2-\norm{u(z_1)-u(z_0)}^2.
\end{equation*}
Take any $w\in\mathbb{R}^n$ and let $z_1=z_0-tw$, $t> 0$. Then the above inequality implies that
\begin{equation*}
\Big\langle T^*\Big(\frac{u(z_0)-u(z_0-tw)}{t}\Big)-w,z_2-z_0\Big\rangle\leq \frac{t\norm{w}^2}2.
\end{equation*}
Letting $t$ tend to zero yields
\begin{equation*}
\langle T^* Du(z_0)w-w,z_2-z_0\rangle\leq 0.
\end{equation*}
As this holds true for any $w\in\mathbb{R}^n$, applying this inequality to $-w$, we infer that the above inequality is an equality, i.e. for all $w\in\mathbb{R}^n$ there is
\begin{equation*}
\langle T^* Du(z_0)w-w,z_2-z_0\rangle= 0.
\end{equation*}
If follows that for all $v\in \mathrm{span}\{z_2-z_0\mid z_2\in\mathcal{S}\}=V$
\begin{equation*}
\langle T^* Du(z_0)w-w,v\rangle= 0,
\end{equation*}
and, consequently, for all such $v$ there is $\langle Q Du(z_0)w- TPw, Tv\rangle=0$. The assertion follows.
\end{proof}

\begin{corollary}\label{col:diff}
Suppose that $\mathcal{S}$ is of dimension $m$. Then $u$ is differentiable in the relative interior of $\mathcal{S}$.
\end{corollary}
\begin{proof}
If the dimension of $\mathcal{S}$ is $m$, then the respective orthogonal projection $Q$ is the identity. The claim follows now by Lemma \ref{lem:diff}.
\end{proof}

\begin{corollary}\label{col:strength}
Let $u\colon\mathbb{R}^n\to\mathbb{R}^m$ be a $1$-Lipschitz map. Let $x_i\in\mathrm{int}\mathcal{S}_i$ belong to the relative interior of leaf $\mathcal{S}_i$ of $u$, for $i=1,2$.  Let $\sigma_i=\mathrm{dist}(\partial\mathcal{S}_i,x_i)$ for $i=1,2$. Then for any $s_1,s_2\in\mathbb{R}^n$ of norm at most one there is
\begin{equation*}
\big|\norm{P_1s_1-P_2s_2}^2-\norm{Du(x_1)s_1-Du(x_2)s_2}^2\big|\leq\frac{\norm{x_1-x_2}^2-\norm{u(x_1)-u(x_2)}^2}{2\sigma_1\sigma_2}.
\end{equation*}
Here $P_i$ denote the orthogonal projection onto the tangent subspace of the leaf $\mathcal{S}_i$ for $i=1,2$. 
Moreover for any $w_1,w_2\in\mathbb{R}^m$ of norm at most one there is
\begin{equation*}
\big|\norm{Q_1w_1-Q_2w_2}^2-\norm{\big(DQ_1u(x_1)\big)^*w_1-\big(DQ_2u(x_2)\big)^*w_2}^2\big|\leq \frac{\norm{x_1-x_2}^2-\norm{u(x_1)-u(x_2)}^2}{2\sigma_1\sigma_2}.
\end{equation*}
Here $Q_i$ denote the orthogonal projection onto the image of $T_i$, for $i=1,2$.
\end{corollary}
\begin{proof}
Formula (\ref{eqn:inequalityimportant}), Lemma \ref{lem:important}, tells us that for any $v_1\in V_1$ and any $v_2\in V_2$ of norm at most one, there is
\begin{equation*}
\big|\norm{v_1-v_2}^2-\norm{T_1v_1-T_2v_2}^2\big|\leq\frac{\norm{x_1-x_2}^2-\norm{u(x_1)-u(x_2)}^2}{2\sigma_1\sigma_2}.
\end{equation*}
Lemma \ref{lem:diff} tells us that $DQ_iu(x_i)=T_iP_i$ for $i=1,2$. Hence the first asserted inequality follows.
Let $v_i=(T_iP_i)^*w_i$ for $w_i\in \mathbb{R}^m$, $i=1,2$, of norm at most one. Then the above formula yields
\begin{equation*}
\big|\norm{(T_1P_1)^*w_1-(T_2P_2)^*w_2}^2-\norm{Q_1w_1-Q_2w_2}^2\big|\leq\frac{\norm{x_1-x_2}^2-\norm{u(x_1)-u(x_2)}^2}{2\sigma_1\sigma_2}.
\end{equation*}
The proof is complete.
\end{proof}

\begin{lemma}\label{lem:boundary}
Let $\mathcal{S}_1, \mathcal{S}_2$ be two distinct leaves of a $1$-Lipschitz map $u\colon \mathbb{R}^n\to\mathbb{R}^m$. Then
\begin{equation*}
\mathcal{S}_1\cap\mathcal{S}_2\subset\partial\mathcal{S}_1\cap\partial\mathcal{S}_2.
\end{equation*}
\end{lemma}
\begin{proof}
We shall first show that there is no point belonging to $ \mathrm{int}\mathcal{S}_1\cap  \mathcal{S}_2$. 
For this, suppose that $x_0\in\mathrm{int}\mathcal{S}_1\cap\mathcal{S}_2$. Let $x_1\in\mathcal{S}_1$ and $x_2\in\mathcal{S}_2$. There exists isometries $T_1$ and $T_2$ on the tangent spaces $V_1$ and $V_2$ of $\mathcal{S}_1$ and $\mathcal{S}_2$ respectively such that
\begin{equation*}
u(x_1)-u(x_0)=T_1(x_1-x_0)\text{ and }u(x_2)-u(x_0)=T_2(x_2-x_0).
\end{equation*}
We may write
\begin{align*}
&\norm{x_1-x_0}^2+\norm{x_2-x_0}^2-2\langle T_1(x_1-x_0),T_2(x_2-x_0)\rangle=\norm{u(x_1)-u(x_2)}^2 \leq\\
&\leq \norm{x_1-x_2}^2=\norm{x_1-x_0}^2+\norm{x_2-x_0}^2-2\langle x_1-x_0,x_2-x_0\rangle.
\end{align*}
Hence
\begin{equation*}
\langle x_1-x_0,x_2-x_0\rangle\leq \langle T_1(x_1-x_0),T_2(x_2-x_0)\rangle.
\end{equation*}
As $x_0\in\mathrm{int}\mathcal{S}_1$ and the inequality holds true for all $x_1\in\mathcal{S}_1$, we actually have equality above for $x_1$ sufficiently close to $x_0$.
It follows that for all $v_1\in V_1$ and $v_2\in V_2$,
\begin{equation}\label{eqn:eqiso}
\langle v_1,v_2\rangle=\langle T_1v_1,T_2v_2\rangle.
\end{equation}
Hence, there exists an isometry that extends both $T_1$ and $T_2$. 
Indeed, define a linear map
\begin{equation*}
S\colon V_1+V_2 \to\mathbb{R}^m
\end{equation*}
by the formula
\begin{equation*}
S(v_1+v_2)=T_1(v_1)+T_2(v_2)\text{ for }v_1\in V_1, v_2\in V_2.
\end{equation*}
We claim that $S$ is a well-defined isometry. Indeed, by (\ref{eqn:eqiso}) and by orthogonality we see that if $v_2\in V_1\cap V_2$, then 
\begin{equation*}
\norm{v_2}^2=\langle T_1v_2,T_2v_2\rangle 
\end{equation*}
which implies, by the equality cases in the Cauchy--Schwarz inequality, that $T_1v_2=T_2v_2$. Thus $S$ is well-defined. It is an isometry, as for $v_1\in V_1$ and $v_2\in V_2$,
\begin{equation*}
\norm{S(v_1+v_2)}^2=\norm{v_1}^2+\norm{v_2}^2+2\langle T_1v_1,T_2v_2\rangle=\norm{v_1+v_2}^2.
\end{equation*}
Moreover, by the definition, $S$ is an extension of both $T_1$ and $T_2$.

Define an affine map $v\colon x_0+V_1 +V_2\to\mathbb{R}^m$ by the formula 
\begin{equation*}
v(x)=S(x-x_0)+u(x_0).
\end{equation*}
Then $v|_{\mathcal{S}_1}=u$ and $v|_{\mathcal{S}_2}=u$.
Choose any points $x\in\mathcal{S}_1$ and $y\in\mathcal{S}_2$. Then 
\begin{equation*}
\norm{u(x)-u(y)}=\norm{v(x)-v(y)}=\norm{S(x-y)}=\norm{x-y}.
\end{equation*}
Thus $u$ is isometric on $\mathcal{S}_1\cup\mathcal{S}_2$. By maximality of leaves, $\mathcal{S}_1=\mathcal{S}_1\cup\mathcal{S}_2=\mathcal{S}_2$, contradicting the distinctness of the two leaves. 
Hence
\begin{equation*}
\mathcal{S}_1\cap\mathcal{S}_2\subset \partial \mathcal{S}_1\cap \mathcal{S}_2.
\end{equation*}
Repeating the above argument with $\mathcal{S}_1$ and $\mathcal{S}_2$ interchanged, we see that 
\begin{equation*}
\mathcal{S}_1\cap\mathcal{S}_2\subset \big(\partial \mathcal{S}_1\cap \mathcal{S}_2\big) \cap \big(\partial \mathcal{S}_2\cap \mathcal{S}_1\big)=\partial \mathcal{S}_1\cap\partial\mathcal{S}_2.
\end{equation*}
\end{proof}

\begin{lemma}\label{lem:notdif}
Let $u\colon\mathbb{R}^n\to\mathbb{R}^m$ be $1$-Lipschitz. If $x_0\in\mathbb{R}^n$ belongs to at least two distinct leaves of $u$, then $u$ is not differentiable at $x_0$.
\end{lemma}
\begin{proof}
Clearly, any zero-dimensional leaf does not intersect any other leaf. Hence, $x_0$ belongs to two distinct leaves $\mathcal{S}_1,\mathcal{S}_2$ of non-zero dimensions. Suppose that $u$ is differentiable at $x_0\in\mathcal{S}_1\cap\mathcal{S}_2$.
Then Lemma \ref{lem:diff} implies that $Q_1u$ is differentiable at $x_0$ with the derivative given by
\begin{equation*}
DQ_1u(x_0)=T_1P_1,
\end{equation*}
where $T_1$ is an isometry such that $u(x)-u(x_0)=T_1(x-x_0)$ for all $x\in\mathcal{S}_1$, $P_1$ is the orthogonal projection onto the tangent space $V_1$ of $\mathcal{S}_1$ and $Q_1$ is the orthogonal projection onto the image of $T_1$. In other words
\begin{equation}\label{eqn:limit}
\lim_{x\to x_0}\frac{Q_1u(x)-Q_1u(x_0)-T_1P_1(x-x_0)}{\norm{x-x_0}}=0.
\end{equation}
For $x\in\mathcal{S}_2$ we may write
\begin{equation*}
u(x)-u(x_0)=T_2(x-x_0)
\end{equation*}
for an isometry $T_2$. If $x\in\mathcal{S}_2$, then
\begin{equation}\label{eqn:difference}
\frac{Q_1u(x)-Q_1u(x_0)-T_1P_1(x-x_0)}{\norm{x-x_0}}=(Q_1T_2-T_1P_1)\bigg(\frac{x-x_0}{\norm{x-x_0}}\bigg).
\end{equation}
For $x\in\mathcal{S}_2$ and $t\in [0,1]$ let
\begin{equation*}
x_t=x_0+t(x-x_0).
\end{equation*}
By convexity of leaves, $x_t\in\mathcal{S}_2$. Observe also that 
\begin{equation}\label{eqn:conv}
\lim_{t\to 0}x_t=x_0.
\end{equation}
It follows by (\ref{eqn:limit}), (\ref{eqn:difference}) and by (\ref{eqn:conv}) that 
\begin{equation*}
Q_1T_2(x-x_0)=T_1P_1(x-x_0)\text{ for all }x\in \mathcal{S}_2.
\end{equation*}
As $V_2=\mathrm{span}\{x-x_0\mid x\in\mathcal{S}_2\}$ is the tangent space of $\mathcal{S}_2$, we infer that 
\begin{equation*}
Q_1T_2v=T_1P_1v\text{ for all }v\in V_2.
\end{equation*}
Hence, for $v_1\in V_1$ and for $v_2\in V_2$
\begin{equation*}
\langle T_1v_1,T_2v_2\rangle=\langle T_1v_1,Q_1T_2v_2\rangle=\langle T_1v_1,T_1P_1v_2\rangle=\langle v_1,v_2\rangle.
\end{equation*}
We continuue the proof as in Lemma \ref{lem:boundary} and arrive at a contradiction that $\mathcal{S}_1=\mathcal{S}_2$.
\end{proof}

\begin{remark}\label{rem:altern}
We may proceed in the first part of the above proof of Lemma \ref{lem:boundary} alternatively. Namely, we may conclude from Lemma \ref{lem:diff} that for at point in $\mathrm{int}\mathcal{S}$ the map $Qu$ is differentiable with $DQu=TP$. Then we proceed as in the proof of Lemma \ref{lem:notdif}.
\end{remark}

\begin{definition}\label{defin:bu}
The set of points belonging to at least two distinct leaves of a $1$-Lipschitz map $u\colon\mathbb{R}^n\to\mathbb{R}^m$ we shall denote by $B(u)$.
\end{definition}

\begin{corollary}\label{col:unique}
For any $1$-Lipschitz function $u\colon\mathbb{R}^n\to\mathbb{R}^m$ the set $B(u)$ is of Lebesgue measure zero.
\end{corollary}
\begin{proof}
Lemma \ref{lem:notdif} implies that $B(u)$ is contained in the set of non-differentiability of $u$. Rademacher's theorem (see e.g. \cite{Federer}) states that the latter is of Lebesgue measure zero.
\end{proof}

\section{Lipschitz change of variables}\label{sec:varia}

Let $u\colon\mathbb{R}^n\to\mathbb{R}^m$ be a $1$-Lipschitz map. The aim of this section is to provide a countable partitioning the union of all $m$-diimensional leaves of $u$ together with a suitable change of variables, which will be useful in the proof of regularity of conditional measures in Section \ref{sec:curv}.

We assume throughout the section that $m\leq n$. Let us recall a lemma taken from \cite[3.2.9]{Federer}.

\begin{lemma}\label{lem:coo}
Let $u\colon\mathbb{R}^n\to\mathbb{R}^m$ be a continuous function. Then the set
\begin{equation*}
\{x\in\mathbb{R}^n\mid u\text{ is differentiable at }x\text{ and }Du(x)\text{ has maximal rank}\}
\end{equation*}
admits a countable Borel cover $(G_i)_{i=1}^{\infty}$  such that for any $i\in\mathbb{N}$ there exists an orthogonal projection $\pi_i\colon\mathbb{R}^n\to\mathbb{R}^{n-m}$ and Lipschitz maps
\begin{equation*}
w_i\colon\mathbb{R}^n\to\mathbb{R}^m\times\mathbb{R}^{n-m}\text{, }v_i\colon\mathbb{R}^m\times\mathbb{R}^{n-m}\to\mathbb{R}^n
\end{equation*}
such that
\begin{equation*}
w_i(x)=(u(x),\pi_i(x))\text{ and }v_i(w_i(x))=x\text{ for all }x\in G_i.
\end{equation*}
\end{lemma}

\begin{lemma}\label{lem:cover}
Let $u\colon\mathbb{R}^n\to\mathbb{R}^m$ be a Lipschitz function. Let $s\in\mathbb{R}^m$ and let
\begin{equation*}
S_s=\{x\in\mathbb{R}^n\mid u(x)=s\}
\end{equation*}
be the level set. Then the set
\begin{equation*}
S_s\cap\{x\in\mathbb{R}^n\mid  u\text{ is differentiable at }x\text{ and }Du(x)\text{ has maximal rank}\}
\end{equation*}
admits a countable Borel covering $(S_s^i)_{i=1}^{\infty}$ of bounded sets such that for all $i\in\mathbb{N}$ there exist Lipschitz functions $w\colon\mathbb{R}^n\to\mathbb{R}^{n-m}$ and $v\colon\mathbb{R}^{n-m}\to\mathbb{R}^n$ satisfying
\begin{equation*}
v(w(x))=x \text{ for all }x\in S_s^i.
\end{equation*}
\end{lemma}
\begin{proof}
We apply Lemma \ref{lem:coo} and obtain a countable covering consisting of Borel sets $G_i$, orthogonal projections $\pi_i\colon\mathbb{R}^n\to\mathbb{R}^{n-m}$ and Lipschitz maps 
\begin{equation*}
w_i\colon\mathbb{R}^n\to\mathbb{R}^m\times\mathbb{R}^{n-m}\text{, }v_i\colon\mathbb{R}^m\times\mathbb{R}^{n-m}\to\mathbb{R}^n
\end{equation*}
such that
\begin{equation*}
w_i(x)=(u(x),\pi_i(x))\text{ and }v_i(w_i(x))=x\text{ for all }x\in G_i.
\end{equation*}
The sets $G_i\cap S_s$ form a countable Borel cover of $S_s$. For any $i\in\mathbb{N}$ define 
\begin{equation*}
w\colon\mathbb{R}^n\to\mathbb{R}^{n-m}\text{ and }
v\colon \mathbb{R}^{n-m}\to\mathbb{R}^n
\end{equation*}
by $w=\pi\circ w_i$, where $\pi\colon \mathbb{R}^m\times\mathbb{R}^{n-m}\to \mathbb{R}^{n-m}$ is the projection on the second variable, and $v(x)=v_i(s,x)$ for $x\in\mathbb{R}^{n-m}$.
Then, if $u(x)=s$, then 
\begin{equation*}
v(w(x))=v_i(s,w(x))=v_i(u(x),\pi_i(x))=x.
\end{equation*}
\end{proof}

\begin{definition}
Pick a countable dense set $S\subset\mathbb{R}^m$. Let $s\in S$. Let $u\colon\mathbb{R}^n\to\mathbb{R}^m$ be a $1$-Lipschitz map. Let $(S_s^i)_{i=1}^{\infty}$ be the Borel cover of Lemma \ref{lem:cover} associated to the level set 
\begin{equation*}
S_s=\{x\in\mathbb{R}^n\mid u(x)=s\}.
\end{equation*}
For each $i,j\in\mathbb{N}$ let the \emph{cluster}
\begin{equation*}
T_{sij}
\end{equation*}
denote the union of all $m$-dimensional leaves $\mathcal{S}$ of $u$ such that there exists $z\in \mathcal{S}\cap S_s^i$ for which  $\mathrm{dist}(z,\partial\mathcal{S})\geq \frac1j$.
Denote by 
\begin{equation*}
T_{sij}^0
\end{equation*}
the union of the relative interiors of all $m$-dimensional leaves $\mathcal{S}$ of $u$ as above.
\end{definition}

\begin{lemma}\label{lem:cluster}
The union of all $m$-dimensional leaves is covered by the clusters 
\begin{equation*}
(T_{sij})_{s\in S,i,j\in\mathbb{N}}.
\end{equation*}
Moreover for each $m$-dimensional leaf $\mathcal{S}$ and each cluster $T_{sij}$ either
\begin{equation*}
\mathrm{int} \mathcal{S}\cap T_{sij}=\emptyset\text{ or }\mathcal{S}\subset T_{sij}.
\end{equation*}
\end{lemma}
\begin{proof}
Let $\mathcal{S}$ be a $m$-dimensional leaf of $u$. Then $u$, if restricted to $\mathcal{S}$, is an isometry onto a convex subset of $\mathbb{R}^m$. Thus, there exists $s\in S\cap \mathrm{int} u(\mathcal{S})$. There exist $i,j\in\mathbb{N}$ and $z\in \mathcal{S}\cap S_s^i$ such that  $\mathrm{dist}(z,\partial\mathcal{S})>1/j$. That is $\mathcal{S}\subset T_{sij}$.

If the interior of some leaf $\mathrm{int}\mathcal{S}$ intersects one of the leaves  comprising the cluster $T_{sij}$, then Lemma \ref{lem:boundary} implies that they are equal and hence $\mathcal{S}\subset T_{sij}$. This completes the proof.
\end{proof}

The lemma below provides the aforementioned change of variables subordinate to the given $1$-Lipschitz map $u\colon\mathbb{R}^n\to\mathbb{R}^m$ and studies it's regularity properties.

\begin{lemma}\label{lem:efge}
Each cluster $T_{sij}\subset\mathbb{R}^n$ admits maps
\begin{equation*}
G\colon T_{sij}^0\to\mathbb{R}^{n-m}\times\mathbb{R}^m
\end{equation*}
and
\begin{equation*}
F\colon G( T_{sij}^0) \to T_{sij}^0
\end{equation*}
such that:
\begin{enumerate}[i)]
\item\label{i:lambda} for each $\lambda>0$, $G$ is a Lipschitz map on the set
\begin{equation*}
T^{\lambda}_{sij}=\Big\{x\in T_{sij}\mid  \mathrm{dist}(x,\partial \mathcal{S}(x))>\lambda\Big\};
\end{equation*}
here $\mathcal{S}(x)$ is the unique leaf of $u$ such that $x\in\mathcal{S}(x)$ and $z\in \mathcal{S}(x)$ is the unique point in $\mathcal{S}(x)$ such that $u(z)=s$,
\item $F$ is Lipschitz on the set $G(T_{sij}^0)$,
\item $F(G(x))=x$ for each $x\in T_{sij}^0$,
\item if a leaf $\mathcal{S}\subset T_{sij}$ intersects $S_s^i$ at a point $z$, then each interior point $x\in \mathrm{int}\mathcal{S}$ of the leaf satisfies
\begin{equation}\label{eqn:giem}
G(x)=(w(z),u(x)-u(z)),
\end{equation}
where $w\colon\mathbb{R}^n\to\mathbb{R}^{n-m}$ is the map from Lemma \ref{lem:cover}.
\end{enumerate}
\end{lemma}
\begin{proof}
Lemma \ref{lem:boundary} shows that the relative interiors of leaves do not intersect any other leaf. Moreover $u$ is an isometry on each leaf. Therefore, every point $x\in T_{sij}^0$ belongs to a unique leaf and each leaf in the cluster $T_{sij}$ intersects the level set $S_s$ in a single point $z\in S_s^i$. It follows that (\ref{eqn:giem}) defines a map 
\begin{equation*}
G\colon T_{sij}^0\to\mathbb{R}^{n-m}\times\mathbb{R}^m,
\end{equation*}
on the cluster $T_{sij}^0$. Let $(a,b)\in G(T_{sij}^0)$ and let $v$ be the map parametrising $S_s^i$ from Lemma \ref{lem:cover}. Then $v(a)\in S_s^i$ belongs to the relative interior of some leaf $\mathcal{S}$ and lies in a distance at least $1/j$ from the relative boundary of the leaf. Define
\begin{equation*}
F(a,b)=v(a)+Du(v(a))^*(b).
\end{equation*}
Let $x\in T_{sij}^0$ belong to a leaf $\mathcal{S}$ that intersects $S_s^i$ at a point $z$. Then $v(w(z))=z$
and there exists an isometry $T$ such that $u(x_1)-u(x_2)=T(x_1-x_2)$ for all $x_1,x_2\in\mathcal{S}$ and $Du(z)=TP$, where $P$ is the orthogonal projection onto the tangent space of $\mathcal{S}$. We infer that
\begin{equation*}
F(G(x))=F(w(z),u(x)-u(z))=z+PT^*T(x-z)=x.
\end{equation*}
We shall now prove that $F$ is Lipschitz on  $G(T_{sij}^0)$. Define
\begin{equation}\label{eqn:lambda}
\Lambda=\big\{a\in\mathbb{R}^{n-m}\mid (a,0)\in G(T_{sij}^0)\big\}.
\end{equation}
We first claim that 
\begin{equation*}
(a,b)\mapsto Du(v(a))^*b
\end{equation*}
is Lipschitz. Recall that $v(a)\in S_s^i$ is in a distance at least $1/j$ from the relative boundary of a leaf $\mathcal{S}$ that contains $v(a)$. Thus, by Corollary \ref{col:strength}, we infer that for points $(a,b),(a',b')\in G(T_{sij}^0)$ there is
\begin{equation*}
\norm{Du(v(a))^*b-Du(v(a'))^*b'}^2
\leq \frac{j^2}{2} \norm{v(a)-v(a')}^2+ \norm{b-b'}^2\leq \frac12C^2j^2\norm{a-a'}^2+\norm{b-b'}^2,
\end{equation*}
where $C$ is the Lipschitz constant of $v$.
It follows immediately that $F$ is Lipschitz on $G(T_{sij}^0)$.

It remains to prove assertion \ref{i:lambda}) of the lemma. Let $\lambda>0$.
Let now $x,x'\in T_{sij}^{\lambda}$ belong to the leaves $\mathcal{S}$ and $\mathcal{S}'$ respectively. By the definition (\ref{eqn:giem}) and by Lipschitzness of $w$ to prove that $G$ is Lipschitz it is enough to show that 
\begin{equation*}
\norm{z-z'}\leq C\norm{x-x'}
\end{equation*}
for some constant $C$. As $u$ is an affine isometry on the leaves we see that
\begin{equation*}
z=x+Du(x)^*(u(z)-u(x))\text{ and }z'=x'+Du(x')^*(u(z')-u(x')).
\end{equation*}
Thus
\begin{equation*}
\norm{z-z'}\leq \norm{x-x'}+\Big\rVert Du(x)^*(u(z)-u(x))-Du(x')^*(u(z')-u(x'))\Big\lVert.
\end{equation*}
Now, by Corollary \ref{col:strength}, taking into account that $u(z)=s=u(z')$, we see that 
\begin{equation*}
\Big\rVert Du(x)^*(u(z)-u(x))-Du(x')^*(u(z')-u(x'))\Big\lVert^2\leq \frac1{2\lambda^2}\norm{x-x'}^2+\norm{u(x)-u(x')}^2.
\end{equation*}
Therefore 
\begin{equation*}
\norm{z-z'}\leq \norm{x-x'}\bigg(1+\sqrt{1+\frac1{2\lambda^2}}\bigg).
\end{equation*}
This concludes the proof that $G$ is Lipschitz on $T_{sij}^{\lambda}$ and completes the proof of the lemma.
\end{proof}

\section{Measurability}\label{sec:measur}

Below $G_{n,k}$ denotes the set of all $k$-dimensional subspaces of $\mathbb{R}^n$. For $V\in G_{n,k}$ we denote by $O_m(V)$ the set of all isometries on $V$ with values in $\mathbb{R}^m$, i.e. the set of all linear maps $T\colon V\to\mathbb{R}^m$ such that 
\begin{equation*}
\norm{T(x)-T(y)}=\norm{x-y}\text{ for all }x,y\in V.
\end{equation*}
By $P_V\colon \mathbb{R}^n\to\mathbb{R}^n$ we denote the orthogonal projection onto $V$. Then $G_{n,k}$ is a compact if equipped with the metric $d$ given by the formula
\begin{equation*}
d(V,V')=\norm{P_V-P_{V'}}\text{, }V,V'\in G_{n,k}.
\end{equation*}
Here $\norm{\cdot}$ denotes the operator norm with respect to the Euclidean norm on $\mathbb{R}^n$.

For a point $x\in\mathbb{R}^n$ and a real number $r>0$ we shall denote by $B(x,r)$ the closed ball centred at $x$ of radius $r$.

\begin{definition}
For $k\in\{1,\dotsc,m\}$ we define $\alpha_k\colon \mathbb{R}^n\to \mathbb{R}\cup\{\infty\}$ by the formula
\begin{equation*}
\alpha_k(x)=\sup\Big\{r\geq 0\mid\exists_{V\in G_{n,k}} \exists_{T\in O_m(V)}\forall_{y\in (x+V)\cap B(x,r)}
 \quad u(x)-u(y)=T(x-y) \Big\}
\end{equation*}
for $x\in\mathbb{R}^n$. 
We define $\alpha_{m+1}\colon\mathbb{R}^n\to\mathbb{R}$ by $\alpha_{m+1}(x)=0$ for all $x\in\mathbb{R}^n$.
\end{definition}

The value of function $\alpha_k(x)$ denotes the greatest radius of a ball such that $u$ is isometric on the intersection of the ball with some $k$-dimensional subspace. 

\begin{lemma}\label{lem:alpha}
For any $k\in\{1,\dotsc,m\}$ the functions $\alpha_k\colon \mathbb{R}^n\to\mathbb{R}\cup\{\infty\}$ are upper semicontinuous.
\end{lemma}
\begin{proof}
Fix $k\in \{1,\dotsc,m\}$. Pick $x_0\in\mathbb{R}^n$ and a sequence $(x_l)_{l=1}^{\infty}$ that converges to $x_0$ such that there exists a limit
\begin{equation*}
\lambda=\lim_{l\to\infty}\alpha_k(x_l).
\end{equation*}
We need to show that $\lambda\leq \alpha_k(x_0)$. Suppose first that $\lambda<\infty$. We may assume that $\alpha_k(x_l)\in\mathbb{R}$ for each $l\in\mathbb{N}$. From the definition of $\alpha_k(x_l)$ it follows that there exist 
\begin{equation*}
V_l\in G_{n,k}\text{ and }T_l\in O_m(V_l)
\end{equation*}
such that for all $y\in (x_l+V_l)\cap B\Big(x_l,\big(1-\frac1l\big)\alpha_k(x_l)\Big)$ we have
\begin{equation*}
u(x_l)-u(y)=T_l(x_l-y).
\end{equation*}
By compactness of $G_{n,k}$ we may assume that the sequence $(V_l)_{l=1}^{\infty}$ is convergent to some $V_0\in G_{n,k}$. Moreover, we may assume that 
\begin{equation*}
(T_lP_{V_l})_{l=1}^{\infty}\text{ converges to }T_0P_{V_0}, 
\end{equation*}
where $T_0\in O_m(V_0)$. Indeed, we may assume that there exists $S_0$ such that $(T_lP_{V_l})_{l=1}^{\infty}$ converges to $S_0$. For $v_0\in V_0$ we have
\begin{equation*}
\norm{v_0}=\lim_{l\to\infty} \norm{P_{V_l}v_0}=\lim_{l\to\infty} \norm{T_lP_{V_l}v_0}=\norm{S_0v_0}.
\end{equation*}
This is to say, $S_0$ is an isometry on $V_0$. As for each $l$ there is $T_lP_{V_l}=T_lP_{V_l}P_{V_l}$, we infer that $S_0=S_0P_{V_0}$. Setting $T_0=S_0P_{V_0}$ proves the claim. 

Choose  now any $v_0\in V_0$ of norm $\norm{v_0}< \lambda$. By the definition of metric on $G_{n,k}$, the sequence $(P_{V_l}v_0)_{l=1}^{\infty}$ converges to $v_0$. Moreover, for sufficiently large $l$, 
\begin{equation*}
x_l+P_{V_l}v_0\in (x_l+V_l)\cap B\Big(x_l,\big(1-\frac1l\big)\alpha_k(x_l)\Big).
\end{equation*}
Thus
\begin{equation*}
u(x_l)-u(x_l+P_{V_l}v_0)=-T_lP_{V_l}v_0.
\end{equation*}
Passing to the limits we obtain that
\begin{equation*}
u(x_0)-u(x_0+v_0)=-T_0v_0.
\end{equation*}
It follows that $\lambda\leq \alpha_k(x_0)$. Thus, the proof is complete provided $\lambda$ is finite. 

Suppose now that $\lambda$ is infinite. Assume again that $\alpha_k(x_l)\in\mathbb{R}$ for each $l\in\mathbb{N}$ and that $(\alpha_k(x_l))_{l=1}^{\infty}$ converges to infinity monotonically. Then there exist $V_l\in G_{n,k}$ and $T_l$ as above, i.e. such that $(V_l)_{l=1}^{\infty}$ converges to $V_0$ and $(T_lP_{V_l})_{l=1}^{\infty}$ converges to $T_0P_{V_0}$, $T_0\in O_m(V_0)$. Taking any $v_0\in V_0$ of norm at most $l\in\mathbb{N}$ we may show that
\begin{equation*}
u(x_0)-u(x_0+v_0)=-T_0v_0.
\end{equation*}
Hence $\alpha_k(x_0)\geq l$ for each $l\in\mathbb{N}$ and thus $\alpha_k(x_0)=\infty$.
\end{proof}

Below we shall denote the unit ball centred at the origin by $B=\{x\in\mathbb{R}^n\mid\norm{x}\leq 1\}$. For $r\geq 0$ we denote by $C_{n,k}(r)$ the set of all $k$-dimensional convex cones $C$ in $\mathbb{R}^n$ such that there exist $c_1,\dotsc,c_k\in C\cap B$ such that the $n\times k$ matrix $D$ with columns $c_1,\dotsc,c_k$ satisfies $\mathrm{det}D^*D\geq r$, i.e. their Gram matrix has determinant at least $r$.
For a cone $C$ we denote by $V_C$ its linear span.

\begin{definition}
For $k\in\{1,\dotsc,m\}$ we define $\beta_k\colon \mathbb{R}^n\to \mathbb{R}$ by the formula
\begin{equation*}
\beta_k(x)=\sup\Big\{r\geq 0\mid\exists_{C\in C_{n,k}(r)}\exists_{T\in O_m(V_{C})}\forall_{y\in (x+C)\cap B(x,r)} \quad
u(x)-u(y)=T(x-y) \Big\}, 
\end{equation*}
where $x\in\mathbb{R}^n$.
Let $\beta_{m+1}\colon \mathbb{R}^n\to\mathbb{R}$ be defined by $\beta_{m+1}(x)=0$ for all $x\in\mathbb{R}^n$.
\end{definition}

The functions $\beta_k(x)$ indicate the maximal radius $r$ such that there is a convex cone $C$, of size in its linear span bounded from below, such that $u$ is isometric on the intersection of the ball centred at $x$ of radius $r$ with the shifted cone $x+C$.

\begin{lemma}\label{lem:upperbeta}
For any $k\in\{1,\dotsc,m\}$ the function $\beta_k\colon \mathbb{R}^n\to\mathbb{R}$ is upper semicontinuous.
\end{lemma}
\begin{proof}
Fix $k\in\{1,\dotsc,m\}$. Pick $x_0\in \mathbb{R}^n$ and a sequence $(x_l)_{l=1}^{\infty}$ that converges to $x_0$ and such that there exists a limit
\begin{equation*}
\lambda=\lim_{l\to\infty}\beta_k(x_l).
\end{equation*}
We need to show that $\lambda\leq \beta_k(x_0)$. Observe that $\lambda<\infty$, as the determinant of a Gram matrix of vectors in the unit ball is bounded above by the volume of the $k$-dimensional unit ball. It follows from the definition of $\beta_k(x_l)$ that there exist 
\begin{equation*}
C_l\in C_{n,k}\Big(\big(1-\frac1l\big)\beta_k(x_l)\Big) \text{ and }T_l\in O_m(V_{C_l})
\end{equation*}
such that for all $y\in (x_l+C_l)\cap B\Big(x_l,\big(1-\frac1l\big)\beta_k(x_l)\Big)$
\begin{equation*}
u(x_l)-u(y)=T_l(x_l-y).
\end{equation*}
For each $l$ pick points $(c_j^l)_{j=1}^k$ in $C_l\cap B$ such that their Gram matrix has determinant at least $(1-1/l)\beta_k(x_l)$. Passing to subsequences, we may assume that the sequences $(c_j^l)_{l=1}^{\infty}$ converge to some points $c_j$, for $j=1,\dotsc,k$, for which the determinant of the Gram matrix is at least $\lambda$. Let $C_0$ be the convex cone in $\mathbb{R}^n$ spanned by $(c_j)_{j=1}^k$, that is
\begin{equation*}
C_0=\Big\{\sum_{j=1}^k\lambda_jc_j\mid\lambda_j\geq 0\text{ for }j=1,\dotsc, k\Big\}.
\end{equation*}
Clearly, $C_0$ has dimension equal to $k$. It follows that $C_0\in C_{n,k}(\lambda)$.

Passing to a subsequence, we may assume that $(V_{C_l})_{l=1}^{\infty}$ converges to some $V_0\in G_{n,k}$. We claim that $V_0=V_{C_0}$.
Choose any $v_0\in V_{C_0}$. Then there exist real numbers $(\lambda_j)_{j=1}^k$ such that
\begin{equation*}
v_0=\sum_{j=1}^k\lambda_jc_j.
\end{equation*}
For $l\in\mathbb{N}$ set $v_l=\sum_{j=1}^k\lambda_jc_j^l$.
Then $(v_l)_{l=1}^{\infty}$ converge to $v_0$ and $v_l\in V_{C_l}$. Hence
\begin{equation*}
v_0=\lim_{l\to\infty} v_l=\lim_{l\to\infty} P_{V_{C_l}}v_l=P_{V_0}v_0.
\end{equation*}
Thus $V_{C_0}\subset V_0$, and the claim follows, as dimension of $V_{C_0}$ is equal to $k$.

As in Lemma \ref{lem:alpha} we show that there exists $T_0\in O_m(V_{C_0})$ such that
\begin{equation*}
(T_lP_{V_{C_l}})_{l=1}^{\infty}\text{ converges to }T_0P_{V_{C_0}}.
\end{equation*}
Take $\epsilon>0$ and choose any $y_0\in (x_0+C_0)\cap B(x_0,(1-\epsilon)\lambda)$. Then there exist $(\lambda_j)_{j=1}^k$ such that
\begin{equation*}
y_0=x_0+\sum_{j=1}^k\lambda_jc_j.
\end{equation*}
Set $y_l=x_l+\sum_{j=1}^k\lambda_jc_j^l$. Then $(y_l)_{l=1}^{\infty}$ converges to $y_0$ and for sufficiently large $l$,
\begin{equation*}
y_l\in(x_l+ C_l)\cap B\Big(x_l,\big(1-\frac1l\big)\beta_k(x_l)\Big).
\end{equation*}
For such $l$ we have $u(x_l)-u(y_l)=T_l(x_l-y_l)$.
It follows that also $u(x_0)-u(y_0)=T_0(x_0-y_0)$.
That is, $\beta_k(x_0)\geq (1-\epsilon)\lambda$ for any $\epsilon>0$. The proof is complete.
\end{proof}

\begin{lemma}\label{lem:beta}
A point $x\in\mathbb{R}^n$ belongs to a leaf $\mathcal{S}$ of $u$ of dimension at least $k$ if and only if $\beta_k(x)>0$. A point $x\in\mathbb{R}^n$ belongs to a leaf $\mathcal{S}$ of $u$ of dimension exactly $k$ if and only if $\beta_k(x)>0$ and $\beta_{k+1}(x)=0$. 
\end{lemma}
\begin{proof}
Suppose that $x_0\in\mathbb{R}^n$ belongs to a leaf $\mathcal{S}$ of $u$ of dimension $l\in\{k,\dotsc,m\}$. Let $V$ denote the tangent space of $\mathcal{S}$. Choose a point $x_1\in\mathrm{int}\mathcal{S}$ and $\epsilon_0>0$ so that the intersection $B(x_1,\epsilon_0)\cap (x_1+V)$ is contained in $\mathcal{S}$. For $\epsilon\in (0,\epsilon_0)$ let 
\begin{equation*}
C=\big\{x\in\mathbb{R}^n\mid x=\lambda (x_2-x_0)\text{ for some }\lambda\geq 0, x_2\in B(x_1,\epsilon)\cap(x_1+ V)\big\}.
\end{equation*}
Then $C$ is a convex cone of dimension $l$ containing the origin. Thus, it contains $k$ linearly independent vectors, which have Gram matrix of non-zero determinant. This is to say, the intersection of $C$ with the linear span of these vectors belongs to $C_{n,k}(\epsilon)$ for $\epsilon>0$ sufficiently small.
Moreover, by convexity of $\mathcal{S}$, $u$ is isometric on the set $(x_0+C)\cap B(x_0,\epsilon)$, if $\epsilon>0$ is sufficiently small. Therefore $\beta_l(x_0)>\epsilon>0$, whenever $\epsilon$ satisfies the two upper bounds.

Conversely, suppose that $\beta_k(x_0)>0$. Then there exist
\begin{equation*}
r>0\text{, a cone }C\in C_{n,k}(r)\text{ and an isometry }T\in O_m(V_{C})
\end{equation*}
such that
\begin{equation*}
u(x_0)-u(y)=T(x-y)\text{ for all }y\in (x_0+C)\cap B(x_0,r).
\end{equation*}
With use of the Kuratowski--Zorn lemma choose a leaf $\mathcal{S}$ of $u$ containing $(x_0+\mathcal{C})\cap B(x_0,\epsilon)$.
Then the dimension of $\mathcal{S}$ is at least $k$. 

The second assertion is a trivial consequence of the first assertion.
\end{proof}

\begin{lemma}\label{lem:interior}
A point $x\in\mathbb{R}^n$ belongs to relative interior of a leaf $\mathcal{S}$ of $u$ of dimension $k$ if and only if $\alpha_k(x)>0$ and $\beta_{k+1}(x)=0$.
\end{lemma}
\begin{proof}
Suppose that $x_0$ belongs to the relative interior of a leaf $\mathcal{S}$ of~$u$ of dimension $k$. By the previous lemma $\beta_k(x_0)>0$ and $\beta_{k+1}(x_0)=0$. 
Let $V$ denote the tangent space of $\mathcal{S}$. Then, as $x_0$ is in the relative interior, there exist $\epsilon>0$,  $T\in O_m(V)$ such that 
\begin{equation*}
u(x_0)-u(y)=T(x_0-y)\text{ for all }y\in (x_0+V)\cap B(x_0,\epsilon).
\end{equation*} 
That is $\alpha_k(x_0)\geq \epsilon>0$. 

Conversely, suppose that $\alpha_k(x_0)>0$ and $\beta_{k+1}(x_0)=0$. 
Then there exist $V\in G_{n,k}$ and $T\in O_m(V)$ such that
\begin{equation*}
u(x_0)-u(y)=T(x_0-y)\text{ for all }y\in (x_0+V)\cap B(x_0,\epsilon).
\end{equation*} 
It follows from the Kuratowski--Zorn lemma that $x_0$ belongs to a leaf $\mathcal{S}$ of $u$ that contains $(x_0+V)\cap B(x_0,\epsilon)$. As $\beta_{k+1}(x_0)=0$, this leaf is of dimension $k$ and thus $x_0$ belongs to the relative interior of $\mathcal{S}$.
\end{proof}

\begin{corollary}\label{col:borel}
Let $k\in\{0,\dotsc,m\}$. Then the union of all leaves of $u$ of dimension $k$ is a Borel set. Moreover, the union of all relative interiors of leaves of $u$ of dimension $k$ is a Borel set and so is the union of all relative boundaries of leaves of $u$ of dimension $k$.
\end{corollary}
\begin{proof}
The proof readily follows by Lemma \ref{lem:alpha}, Lemma \ref{lem:upperbeta} and  Lemma \ref{lem:beta}, Lemma \ref{lem:interior}.
\end{proof}

Note that whenever $\alpha_k(x)>0$ and $\beta_{k+1}(x)=0$, then Lemma \ref{lem:interior} tells us that $x$ belongs to the relative interior of a leaf of $u$. This leaf in unique, by Lemma \ref{lem:boundary}. We shall denote it by $\mathcal{S}(x)$.

Below we adapt a convention that $\inf\emptyset=\infty$. 

\begin{definition}\label{def:gamma}
Let $k\in \{0,\dotsc,m\}$. We define $\gamma_k\colon \mathbb{R}^n\times\mathbb{R}^m\to \mathbb{R}\cup\{\infty\}$ by the formula
\begin{equation*}
\gamma_k(x,y)=\inf\Big\{t>0\mid y\in t\big(u(\mathcal{S}(x))-u(x)\big)\Big\}
\end{equation*}
for $x\in\mathbb{R}^n$ such that $\alpha_k(x)>0$ and $\beta_{k+1}(x)=0$
and
\begin{equation*}
\gamma_k(x,y)=\infty
\end{equation*}
otherwise.
\end{definition}

\begin{lemma}
For any $k\in \{0,\dotsc,m\}$ the function $\gamma_k$ is Borel measurable.
\end{lemma}
\begin{proof}
As $\alpha_k$ and $\beta_{k+1}$ are Borel measurable, it is enough to show that the function $\gamma_k$ is Borel measurable on 
\begin{equation*}
A_k=\big\{(x,y)\in\mathbb{R}^n\times\mathbb{R}^m\mid\alpha_k(x)>0\text{ and }\beta_{k+1}(x)=0\big\}.
\end{equation*}
Observe that $\gamma_k$ is a limit, as $\rho$ converges to infinity, of functions
\begin{equation*}
\gamma_{k,\rho}(x,y)=\inf\Big\{t>0\mid y\in t\big(u(\mathcal{S}(x))-u(x)\big),\norm{y}\leq t\rho\Big\}.
\end{equation*}
We claim that $\gamma_{k,\rho}$ is lower semicontinuous on $A_k$. This will yield the asserted measurability.

Indeed, let $(x_l,y_l)_{l=1}^{\infty}$ be a sequence in $A_k$ such that there exists $(x_0,y_0)\in A_k$ with
\begin{equation*}
(x_0,y_0)=\lim_{l\to\infty} (x_l,y_l) \text{ and such that there exists }\lim_{l\to\infty}\gamma_{k,\rho}(x_l,y_l)=\lambda.
\end{equation*}
We shall show that $\gamma_{k,\rho}(x_0,y_0)\leq\lambda$. If $\lambda=\infty$, then there is nothing to prove. Otherwise, there exist sequences $(z_l)_{l=1}^{\infty}$ in $\mathbb{R}^n$ and $(t_l)_{l=1}^{\infty}$ in $\mathbb{R}$ such that 
\begin{equation}\label{eqn:equality3}
y_l=t_l\big(u(z_l)-u(x_l)\big)\text{ and }\norm{y_l}\leq t_l \rho\text{, where }z_l\in \mathcal{S}(x_l)\text{ and }0<t_l<\gamma_k(x_l,y_l)+1/l.
\end{equation}
 Observe that 
\begin{equation*}
\norm{z_l-x_l}=\norm{u(z_l)-u(x_l)}=\frac{\norm{y_l}}{t_l}\leq  \rho
\end{equation*}
Thus, passing to a subsequence, we may assume that $(z_l)_{l=1}^{\infty}$ converges to some $z_0\in\mathcal{S}(x_0)$ and that $(t_l)_{l=1}^{\infty}$ converges to some $t_0\geq 0$. Taking limits in (\ref{eqn:equality3}) we see that
\begin{equation*}
y_0=t_0\big(u(z_0)-u(x_0)\big)\text{ with }z_0\in \mathcal{S}(x_0)\text{ and }0\leq t_0\leq \lambda.
\end{equation*}
Hence
\begin{equation*}
y_0\in t_0 \big(u(\mathcal{S}(x_0))-u(x_0)\big)\text{ and }\norm{y_0}\leq t_0\rho.
\end{equation*}
This is to say, $\gamma_{k,\rho}(x_0,y_0)\leq t_0\leq\lambda$.
The proof is complete.
\end{proof}

\begin{definition}
For a convex set $K\subset\mathbb{R}^m$, such that $0\in \mathrm{int}K$, we define its \emph{Minkowski functional} $\norm{\cdot}_K\colon\mathbb{R}^m\to\mathbb{R}\cup\{\infty\}$ by the formula
\begin{equation*}
\norm{y}_K=\inf\big\{t>0\mid y\in tK\big\}.
\end{equation*}
\end{definition}

The following proposition can be found e.g. in \cite[Theorem 5.3.3]{Narici}.

\begin{proposition}\label{pro:minkowski}
Let $K\subset \mathbb{R}^m$ be a closed, convex set that contains the origin in its interior. A point $y\in\mathbb{R}^m$ belongs to the interior of $K$ if and only if $\norm{y}_K<1$.

Moreover, a point $y\in \mathbb{R}^m$ belongs to the boundary of $K$ if and only if $\norm{y}_K=1$.
\end{proposition}

%
%
%

\begin{lemma}\label{lem:minkowski}
If $x\in\mathbb{R}^n$ belongs to relative interior of a leaf $\mathcal{S}$ of $u$ of dimension $k$, then $\gamma_k(x,\cdot)$ is the Minkowski functional of the closed, convex set $u(\mathcal{S})-u(x)$.
If $x\in\mathbb{R}^n$ does not belong to relative interior of any leaf of dimension $k$, then 
\begin{equation*}
\gamma_k(x,\cdot)=\infty.
\end{equation*}
\end{lemma}
\begin{proof}
Suppose that $x\in\mathbb{R}^n$ does not belong to relative interior of a leaf of $u$ of dimension $k$. Then Lemma \ref{lem:interior} and Definition \ref{def:gamma} tells us that $\gamma_k(x,\cdot)=\infty$. 

Let now $x\in\mathrm{int}\mathcal{S}$, where $\mathcal{S}$ is a $k$-dimensional leaf. By Lemma \ref{lem:boundary}, such leaf $\mathcal{S}$ is unique. The assertion of the lemma follows readily from the definitions.
\end{proof}

\begin{definition}
Let $k\in \{0,\dotsc,m\}$. We shall denote by $T_k$ the union of all $k$-dimensional leaves of $u$, by $\mathrm{int}T_k$ the union of all relative interiors of all $k$-dimensional leaves of $u$ and by $\partial{T}_k$ the union of all relative boundaries of all $k$-dimensional leaves of $u$.
\end{definition}

Below we shall denote by $\lambda$ the Lebesgue measure. The space on which $\lambda$ is considered will be clear from the context.

\begin{lemma}\label{lem:measurable}
For each $s\in S$ and each $i,j\in\mathbb{N}$ the cluster $T_{sij}^0$ and its image $G(T_{sij}^0)$ are Borel sets. Moreover $\partial T_m$ is a Borel set of Lebesgue measure zero.
\end{lemma}
\begin{proof}
Fix $s\in S$ and $i,j\in\mathbb{N}$. Recall the Borel set $S_s^i\subset\mathbb{R}^n$ and Lipschitz mapping $w\colon \mathbb{R}^n\to\mathbb{R}^{n-m}$ from Lemma \ref{lem:cover}. Since $w$ is injective on $S_s^i$ it follows from \cite[2.2.10]{Federer} that $w(S_s^i)$ is a Borel subset of $\mathbb{R}^{n-m}$. Moreover, the set $\Lambda$, defined in (\ref{eqn:lambda}), is given by
\begin{equation}\label{eqn:lam}
\Lambda=\Big\{a\in w(S_s^i)\mid \alpha_m(w^{-1}(a))>1/j\Big\}
\end{equation}
as follows by the definition (\ref{eqn:giem}) and Lemma \ref{lem:cover}. 
Clearly, $\Lambda$ is a Borel set. 
Definition of the cluster $T_{sij}^0$ implies that
\begin{equation*}
G(T_{sij}^0)=\bigg\{(a,b)\in\mathbb{R}^{n-m}\times\mathbb{R}^m\mid a\in \Lambda, b\in u\Big(\mathrm{int}\mathcal{S}\big(v(a)\big)\Big)-u\big(v(a)\big)\bigg\}.
\end{equation*}
Here $\mathcal{S}(v(a))$ is the unique $m$-dimensional leaf of $u$ containing $v(a)$.
Observe that Proposition \ref{pro:minkowski} and Lemma \ref{lem:minkowski} tells us that if $a\in\Lambda$, then $b$ belongs to the interior of
\begin{equation*}
u(\mathcal{S}(v(a)))-u(v(a))\text{ if and only if }
\gamma_m(v(a),b)<1.
\end{equation*}
This is to say,
\begin{equation}\label{eqn:image}
G(T_{sij}^0)=\Big\{(a,b)\in\mathbb{R}^{n-m}\times\mathbb{R}^m\mid a\in \Lambda,\gamma_m(v(a),b)<1\Big\}.
\end{equation}
As $\gamma_m$ is Borel measurable, it follows that $G(T_{sij}^0)$ is a Borel set.

Lemma \ref{lem:efge} shows that $F$, the inverse of $G$ on its image, is well-defined, injective and Lipschitz on $G(T_{sij}^0)$. 
Moreover
\begin{equation*}
T_{sij}^0=F(G(T_{sij}^0)).
\end{equation*}
Using \cite[2.2.10]{Federer}, we see that $T_{sij}^0$ is a Borel set. 

We shall show that $\partial T_m$ has Lebesque measure zero. Recall that Corollary \ref{col:borel} tells us that $\partial T_m$ is a Borel set. 
Consider the set
\begin{equation*}
B=\Big\{(a,b)\in\mathbb{R}^{n-m}\times\mathbb{R}^m\mid a\in\Lambda,\gamma_m(v(a),b)=1\Big\}.
\end{equation*}
By Fubini's theorem, $\lambda(B)=0$, as boundaries of convex sets have Lebesgue measure zero. 

Recall that $F$ is a Lipschitz map on $G(T_{sij}^0)$. Using the Kirszbraun theorem (see e.g \cite{Kirszbraun, Schoenberg}) we extend $F$, defined on $G(T_{sij}^0)$, to a Lipschitz map $\tilde{F}$ on $\mathbb{R}^{n-m}\times\mathbb{R}^m$.
We claim that for any such extension
\begin{equation}\label{eqn:exten}
\tilde{F}(B)\supset\partial T_m.
\end{equation}
Indeed, let $x\in \partial T_m$. 
There exists a leaf $\mathcal{S}\subset T_{sij}$ of $u$ and a sequence $(x_l)_{l=1}^{\infty}$ in $\mathrm{int}S$ that converges to $x$. Let $\tilde{G}$ be a Lipschitz extension of $G$ to $\mathbb{R}^n$. The sequence $(G(x_l))_{l=1}^{\infty}$ converges to $\tilde{G}(x)=(a,b) \in\mathbb{R}^{n-m}\times\mathbb{R}^m$. We claim that $(a,b)\in B$. This follows by the continuity of $\gamma_m$ in the second variable. Now, $x=\tilde{F}(\tilde{G}(x))\in \tilde{F}(B)$ and (\ref{eqn:exten}) is proven.

Therefore we can use $\lambda(B)=0$ and the fact that images under Lipschitz maps of sets of Lebesgue measure zero have Lebesgue measure zero (see \cite[3.2.3]{Federer}), to infer that $\lambda(\partial T_m\cap T_{sij} )=0$
and hence $\partial T_m\cap T_{sij}$ is Lebesgue measurable. By Lemma \ref{lem:cluster} the sets $T_{sij}$ form a countable cover of $\partial T_m$. It follows that $\lambda(\partial T_m)=0$. This concludes the proof.
\end{proof}

\begin{corollary}
For any $s\in S$, $i,j\in\mathbb{N}$, the set $T_{sij}$ is Lebesgue measurable.
\end{corollary}
\begin{proof}
$T_{sij}$ is a union of a Borel set $T_{sij}^0$ and a set $\partial T_m\cap T_{sij}$ of Lebesgue measure zero.
\end{proof}

\begin{remark}\label{rem:cluster}
The clusters $T_{sij}$ may be taken to be disjoint. Indeed, let $(T_k)_{k=1}^{\infty}$ be a renumbering of the set of clusters. Set for $l\in\mathbb{N}$
\begin{equation*}
T_l'=T_l\setminus \bigcup_{j=1}^{l-1}T_j
\end{equation*} 
and 
\begin{equation*}
\mathrm{int}T_l'=\mathrm{int}T_l\setminus \bigcup_{j=1}^{l-1}\mathrm{int}T_j.
\end{equation*} 
Note that the structure of the clusters $T'_{sij}$ remains the same. For each $T_{sij}$ there exists a Borel subset $S_{sij}=T_{sij}\cap S_s^i$ of $S_s^i\subset\mathbb{R}^n$ on which there are Lipschitz maps 
\begin{equation*}
w\colon\mathbb{R}^n\to\mathbb{R}^{n-m}\text{ and }v\colon\mathbb{R}^{n-m}\to\mathbb{R}^n\end{equation*}
such that
\begin{equation*}
v(w(x))=x\text{ for all } x\in S_{sij}
\end{equation*}
Indeed, the new cluster is a subset of the old one, so the former maps suffice. From the modification procedure it follows also that Lemma \ref{lem:cluster} still holds true. Moreover, the leaf $\mathcal{S}$ corresponding to a point $z\in S_{sij}$ satisfies
\begin{equation*}
\mathrm{dist}(z,\partial\mathcal{S})>1/j.
\end{equation*}
Also the assertions of Lemma \ref{lem:efge} hold true with the old maps and so does the assertions of Lemma \ref{lem:measurable}, as follows from the modification procedure.
\end{remark}

\section{Disintegration with respect to partition}\label{sec:disin}

Let $u\colon\mathbb{R}^n\to\mathbb{R}^m$ be a $1$-Lipschitz map with respect to the Euclidean norms. In the previous sections we have associated to $u$ a partitioning of $\mathbb{R}^n$, up to a set of Lebesgue measure zero, into maximal sets $\mathcal{S}$ on which $u$ is an isometry. It was conjectured by Klartag in \cite[Chapter 6]{Klartag} that given a measure $\mu$, such that $(\mathbb{R}^n,\norm{\cdot},\mu)$ is a weighted Riemannian manifold satisfying the curvature-dimension condition $CD(\kappa,N)$ (see Definition \ref{defin:curv}), then $\mu$ may be decomposed into a mixture of measures $\mu_{\mathcal{S}}$, each supported on a leaf $\mathcal{S}$ of $u$, such that $(\mathrm{int}\mathcal{S},\norm{\cdot},\mu_{\mathcal{S}})$ is a weighted Riemannian manifold that satisfies $CD(\kappa,N)$.

Below we denote by $CC(\mathbb{R}^n)$ the space of closed, convex, non-empty subsets of $\mathbb{R}^n$. It is a closed subspace of $CL(\mathbb{R}^n)$ -- the space of closed non-empty subsets of $\mathbb{R}^n$ equipped with the Wijsman topology (see \cite{Wijsman}). The Wijsman topology is the weakest topology such that for any $x\in\mathbb{R}^n$ function
\begin{equation*}
A\mapsto \mathrm{dist}(x,A)
\end{equation*}
is continuous. By a result of Beer (see \cite{Beer2}), the space $CL(\mathbb{R}^n)$, equipped with this topology, is Polish. Hence so is $CC(\mathbb{R}^n)$.

Let us recall that $B(u)$ denotes the set of points in $\mathbb{R}^n$ that belong to at least two distinct leaves of $u$. By Corollary \ref{col:unique}, it is contained in the Borel set $N(u)$ of points at which $u$ is not differentiable. The latter is of Lebesgue measure zero. We define a map
\begin{equation*}
\mathcal{S}\colon\mathbb{R}^n\to CC(\mathbb{R}^n)
\end{equation*}
in such a way that for $x\in\mathbb{R}^n\setminus N(u)$ the set $\mathcal{S}(x)$ is the unique leaf of $u$ containing $x$ and for $x\in N(u)$ we put $\mathcal{S}(x)=\{x\}$.

The aim of this section is to prove the following disintegration theorem, which is a step towards the conjecture. 

\begin{theorem}\label{thm:dis}
Let $u\colon\mathbb{R}^n\to\mathbb{R}^m$ be a $1$-Lipschitz map with respect to the Euclidean norms. Then there exists a Borel measure $\nu$ on $CC(\mathbb{R}^n)$, supported on the set of leaves of $u$, and Borel measures $\lambda_{\mathcal{S}}$ such that 
\begin{enumerate}[i)]
\item for every Borel set $A\subset\mathbb{R}^n$ the function $\mathcal{S}\mapsto \lambda_{\mathcal{S}}(A)$ is Borel measurable,
\item for $\nu$-almost every leaf $\mathcal{S}$ the measure $\lambda_{\mathcal{S}}$ is concentrated on $\mathcal{S}$,
\item  for every Borel set $A\subset\mathbb{R}^n$
\begin{equation*}
\lambda(A)=\int_{CC(\mathbb{R}^n)} \lambda_{\mathcal{S}}(A)d\nu(\mathcal{S}).
\end{equation*}
\end{enumerate}
\end{theorem}

Let $X$ be a measurable space. In \cite{Beer3} it is proven that a map $f\colon X\to CL(\mathbb{R}^n)$ is measurable if and only if it is measurable as a multifunction. The latter is defined by the condition that for any open set $U\subset\mathbb{R}^n$ the set
\begin{equation*}
\{x\in X\mid f(x)\cap U\neq \emptyset\}
\end{equation*}
is measurable in $X$.
 
Let us recall a theorem that follows readily from \cite[Example 10.4.11, Definition 10.4.1]{Bogachev2}. 

\begin{theorem}\label{thm:disinteg}
Let $X,Y$ be two Polish spaces. Let $\pi\colon X\to Y$ be a Borel map and let $\mu$ be a non-negative finite Borel measure on $X$. Let $\nu$ be the push-forward of measure $\mu$ via $\pi$. Then there exist Borel measures $(\mu_y)_{y\in Y}$ on $X$ such that
\begin{enumerate}[i)]
\item  for every Borel set $B\subset X$ the function $y\mapsto\mu_y(B)$ is Borel measurable,
\item for $\nu$-almost every $y\in \pi(X)$ the measure $\mu_y$ is concentrated on $\pi^{-1}(y)$,
\item for every Borel sets $B\subset X$ and $E\subset Y$ there is
\begin{equation*}
\mu(B\cap \pi^{-1}(E))=\int_E \mu_y(B) d\nu(y).
\end{equation*}
\end{enumerate}
\end{theorem}

\begin{proof}[Proof of Theorem \ref{thm:dis}]
We have a well-defined map $\mathcal{S}\colon \mathbb{R}^n\to CC(\mathbb{R}^n)$ that assigns to any $x\in \mathbb{R}^n\setminus N(u)$ a unique leaf $\mathcal{S}(x)$ that contains $x$ and for $x\in N(u)$ we set $\mathcal{S}(x)=\{x\}$.
We would like to prove that $\mathcal{S}\colon \mathbb{R}^n\to CC(\mathbb{R}^n)$ is Borel measurable with respect to the Wijsman topology on $CC(\mathbb{R}^n)$, which is equivalent to its measurablity as a multifunction.
Note that for any compact set $K\subset\mathbb{R}^n$ the set $A_K=\{x\in \mathbb{R}^n\mid \mathcal{S}(x)\cap K\neq \emptyset\}$ is equal to
\begin{equation*}
 \Big\{x\in \mathbb{R}^n\setminus  (K\cup N(u))\mid \sup\Big\{\frac{\norm{u(x)-u(y)}}{\norm{x-y}}\mid y\in K\Big\}=1\Big\}\cup K.
\end{equation*}
Observe that the function
\begin{equation*}
x\mapsto  \sup\Big\{\frac{\norm{u(x)-u(y)}}{\norm{x-y}}\mid y\in K\Big\}
\end{equation*}
is lower semicontinuous. Hence $A_K$ is a Borel set. As any open set $U\subset\mathbb{R}^n$ is a countable union of compact sets, it follows that the map $\mathcal{S}$ is Borel measurable.

Recall that $CC(\mathbb{R}^n)$ and $\mathbb{R}^n$ are Polish spaces.
We partition $\mathbb{R}^n$ into countably many closed sets of finite Lebesgue measure such that the measure of the set of points that belong to at least two elements of this partition is zero. To each element of the partition we apply Theorem \ref{thm:disinteg}. Summing up the resulting conditional measures, and taking into account that the set $N(u)$ has Lebesgue measure zero, we obtain the desired disintegration.
\end{proof}

\section{Curvature-dimension condition}\label{sec:curv}

Suppose that we are given a measure $\mu$ on $\mathbb{R}^n$ such that $(\mathbb{R}^n,\norm{\cdot},\mu)$ is a weighted Riemannian manifold satisfying $CD(\kappa, N)$. We shall investigate the behaviour of the conditional measures of $\mu$, see Section \ref{sec:disin}, with respect to the partition introduced in Section \ref{sec:partition}. We shall concentrate on the leaves of maximal dimension.

Let us recall the notion of the curvature-dimension condition $CD(\kappa,n)$. We shall say that an $n$-dimensional Riemannian manifold $\mathcal{M}$ satisfies the $CD(\kappa,n)$ condition provided that the Ricci tensor $Ric_M$ is bounded below by the Riemannian metric tensor $g$, i.e.
\begin{equation*}
Ric_{\mathcal{M}}(p)(v,v)\geq\kappa g(p)(v,v)\text{ for any }p\in \mathcal{M}\text{ and any }v\in T_p\mathcal{M}.
\end{equation*}
We shall study weighted Riemannian manifolds, which are triples $(\mathcal{M},d,\mu)$, where $d$ is the Riemannian metric on $\mathcal{M}$ and $\mu$ is a measure on $\mathcal{M}$ with smooth positive density $e^{-\rho}$ with respect to the Riemannian volume. The generalised Ricci tensor of the weighted Riemannian manifold is defined by the formula
\begin{equation*}
Ric_{\mu}=Ric_\mathcal{M}+D^2\rho,
\end{equation*}
where $D^2\rho$ is the Hessian of smooth function $\rho$. The generalised Ricci tensor -- or the $N$-Bakry-\'Emery tensor -- with parameter $N\in (-\infty,1)\cup [n,\infty]$ is defined by the formula
\begin{equation*}
Ric_{\mu,N}(v,v)=\begin{cases}
    Ric_{\mu}(v,v)-\frac{D\rho(v)^2}{N-n},& \text{if } N>n\\
    Ric_{\mu}(v,v)&\text{if }N=\infty\\
    Ric_{\mathcal{M}}(v,v)    &\text{if }N=n\text{ and }\rho\text{ is constant.}
\end{cases}
\end{equation*}
Note that if $N=n$, then $\rho$ is required to be a constant function.

\begin{definition}\label{defin:curv}
For $\kappa\in\mathbb{R}$ and $N\in (-\infty,1)\cup [n,\infty]$ we say that $(\mathcal{M},d,\mu)$ satisfies the curvature-dimension condition $CD(\kappa,N)$ if 
\begin{equation*}
Ric_{\mu,N}(p)(v,v)\geq \kappa g(p)(v,v)\text{ for all }p\in \mathcal{M}\text{ and all }v\in T_p\mathcal{M}.
\end{equation*}
\end{definition}

We refer the reader to \cite{Bakry1}, \cite{Bakry},  \cite{Ledoux} and to \cite{Ambrosio4}, \cite{Sturm3} \cite{Sturm1}, \cite{Sturm2}, \cite{Villani1}  for background on the curvature-dimension condition. 
In all cases we consider in this paper it will always hold that $Ric_{\mathcal{M}}=0$. 

The aim of the section is to prove the following theorem which partially resolves the conjecture of Klartag \cite[Chapter 6]{Klartag} in the affirmative. In particular, if a measure $\mu$ is concentrated on leaves of $u$ of dimension $m$, then the conjecture holds true for $\mu$ and $u$.

Let us recall that $T_m$ denotes the union of leaves of dimension $m$. This is a Borel set by Corollary \ref{col:borel}.

We shall denote by $T^m\subset CC(\mathbb{R}^n)$ the set of leaves of $u$ of dimension $m$.

Below we present a generalisation of Theorem \ref{thm:localisation}.

\begin{theorem}\label{thm:discd}
Let $m\leq n$. Let $N\in (-\infty,1)\cup [n,\infty]$ andl let $\kappa\in\mathbb{R}$. Let $u\colon\mathbb{R}^n\to\mathbb{R}^m$ be a $1$-Lipschitz map with respect to the Euclidean norms. Let $\mu$ be a Borel measure on $\mathbb{R}^n$ such that $(\mathbb{R}^n,\norm{\cdot},\mu)$ satisfies the curvature-dimension condition $CD(\kappa,N)$. 
Then there exists a Borel measure $\nu$ on $CC(\mathbb{R}^n)$, supported on the set $T^m$ of leaves of dimension $m$, and for each leaf $\mathcal{S}$ of $u$ of dimension $m$, there exists a Borel measure $\mu_{\mathcal{S}}$ such that:
\begin{enumerate}[i)]
\item\label{i:measur} for every Borel set $B\subset T_m$ the function $\mathcal{S}\mapsto \mu_{\mathcal{S}}(B)$ is $\nu$-measurable,
\item\label{i:concent} for $\nu$-almost every leaf $\mathcal{S}$ the measure $\mu_{\mathcal{S}}$ is concentrated on $\mathrm{int}\mathcal{S}$
\item\label{i:curv} for $\nu$-almost every leaf $\mathcal{S}$ the space $(\mathrm{int}\mathcal{S},\norm{\cdot},\mu_{\mathcal{S}})$ satisfies the $CD(\kappa,N)$ condition,
\item\label{i:disinteg} for every Borel set $A\subset T_m$ there is
\begin{equation*}
\mu(A)=\int_{T^m} \mu_{\mathcal{S}}(A)d\nu(\mathcal{S}).
\end{equation*} 
\end{enumerate}
\end{theorem}

Let us note that in \cite[Lemma 3]{Ciosmak1} it is proven that the set of trivial leaves of a $1$-Lipschitz map is of Lebesgue measure zero. Therefore, the above theorem is indeed a generalisation of Theorem \ref{thm:localisation}.

In what follows, we shall use the notation from Section \ref{sec:measur}. Observe that it suffices to prove the theorem under the assumption that $\mu$ is concentrated on a single cluster $T_{sij}$, $s\in S$ and $i,j\in\mathbb{N}$, of leaves of $u$; see Lemma \ref{lem:measurable} and Remark \ref{rem:cluster}. Recall the definitions of maps $F$ and $G$ (see Lemma \ref{lem:efge}) and a map $v$ (see Lemma \ref{lem:cover}). Below $\mathcal{H}_m$ is the $m$-dimensional Hausdorff measure on $\mathbb{R}^n$. 

\begin{lemma}\label{lem:density}
Let $m\leq n$ and let $u\colon\mathbb{R}^n\to\mathbb{R}^m$ be a $1$-Lipschitz map. Fix $s\in S$, $i,j\in\mathbb{N}$. Then for any Borel set $A\subset T_{sij}$ there is
\begin{equation*}
\lambda(A)=\int_{\Lambda}\Big(\int_{\mathrm{int}\mathcal{S}(v(a))}\mathbf{1}_A J_nF\circ G d\mathcal{H}_m\Big) d\lambda(a),
\end{equation*}
where $J_nF$ denotes the $n$-dimensional Jacobian of $F$ and 
\begin{equation*}
\Lambda=\big\{a\in\mathbb{R}^{n-m}\mid (a,0)\in G(T_{sij}^0)\big\}.
\end{equation*}
Moreover, the map
\begin{equation*}
\Lambda\ni a\mapsto \int_{\mathrm{int}\mathcal{S}(v(a))}\mathbf{1}_A J_nF\circ G d\mathcal{H}_m\in\mathbb{R}
\end{equation*}
is $\lambda$-measurable.
\end{lemma}
\begin{proof}
By Lemma \ref{lem:efge}, the map $F$ is a bijection of $G(T_{sij}^0)$ and of $T_{sij}^0$. As $F$ is Lipschitz on $G(T_{sij}^0)$ we may apply the area formula \cite[3.2.5]{Federer} to infer that for any measurable, non-negative $\phi\colon\mathbb{R}^n\to\mathbb{R}$ 
\begin{equation}\label{eqn:area}
\int_{G(T_{sij}^0)}\phi\circ F J_nF d\lambda=\int_{T_{sij}^0}\phi d\lambda.
\end{equation}
Let
\begin{equation*}
f=J_nF\mathbf{1}_{G(T_{sij}^0)}.
\end{equation*}
Observe that $f$ is non-negative and Borel measurable as $G(T_{sij}^0)$ is a Borel set by Lemma \ref{lem:efge} and the fact that images of Borel sets via Lipschitz maps are Borel.

By Tonelli's theorem, the functions $f(a,\cdot)$ are measurable for almost every $a\in\mathbb{R}^{n-m}$ and we have
\begin{equation}\label{eqn:fubini}
\int_{\mathbb{R}^{n-m}\times\mathbb{R}^m}\phi\circ Ffd\lambda=\int_{\mathbb{R}^{n-m}}\int_{\mathbb{R}^m}\phi (F(a,b))f(a,b)d\lambda(b)d\lambda(a).
\end{equation}
Observe now that $(a,b)\in G(T_{sij}^0)$ if and only if for some $a\in\Lambda$
\begin{equation*}
a=w(v(a))\text{ and }b=u(x)-u(v(a)).
\end{equation*}
Note that $F$ on $G(\mathrm{int}\mathcal{S}(v(a)))$ is an isometry. Therefore by a linear change of variables
\begin{equation*}
\int_{G(\mathrm{int}\mathcal{S}(v(a)))}\phi(F(a,b))f(a,b)d\lambda(b)=\int_{\mathrm{int}\mathcal{S}(v(a))}\phi f\circ G d\mathcal{H}_m.
\end{equation*}
Tonelli's theorem implies that the map
\begin{equation*}
\Lambda\ni a\mapsto \int_{\mathrm{int}\mathcal{S}(v(a))}\phi f\circ G d\mathcal{H}_m
\end{equation*}
is measurable. Moreover, by (\ref{eqn:area}) and by (\ref{eqn:fubini}), for any non-negative function $\phi$ we have
\begin{equation*}
\int_{T_{sij}^0}\phi d\lambda=\int_{\Lambda}\Big(\int_{\mathrm{int}\mathcal{S}(v(a))}\phi f\circ G d\mathcal{H}_m\Big) d\lambda(a).
\end{equation*}
By the fact that $\partial T_m$ has Lebesgue measure zero (see Lemma \ref{lem:measurable}), we see that 
\begin{equation*}
\int_{T_{sij}}\phi d\lambda=\int_{\Lambda}\Big(\int_{\mathrm{int}\mathcal{S}(v(a))}\phi f\circ G d\mathcal{H}_m\Big) d\lambda(a).
\end{equation*}
The proof is complete.
\end{proof}

Let us recall a lemma from \cite{Klartag} that we shall need in what follows.

\begin{lemma}\label{lem:tri}
Let $a,b\in\mathbb{R}$, $b>0$ and $a\notin [-b,0]$. Then
\begin{equation*}
\frac{x^2}{a}+\frac{y^2}{b}\geq\frac{(x-y)^2}{a+b}
\end{equation*}
for all $x,y\in\mathbb{R}$.
\end{lemma}
\begin{proof}
We use the inequality
\begin{equation*}
\frac{\abs{a}}{\abs{b}}x^2\pm 2xy+\frac{\abs{b}}{\abs{a}}y^2\geq 0.
\end{equation*}
From this we see that
\begin{equation*}
\frac{x^2}{a}+\frac{y^2}{b}-\frac{(x-y)^2}{a+b}= \frac1{a+b}\Big(\frac{b}{a}x^2+2xy+\frac{a}{b}y^2\Big)\geq 0
\end{equation*}
whenever $b>0$ and $a\notin [-b,0]$.
\end{proof}

Let us also recall formulae for differentiation of matrices. If $R(t)= \log\abs{\det A(t)}$ and $A$ is differentiable in $t\in\mathbb{R}$, then
\begin{equation}\label{eqn:formula}
\frac{dR}{dt}(s)=\mathrm{tr}\Big(A(s)^{-1}\frac{ dA}{dt}(s)\Big).
\end{equation}
Moreover
\begin{equation}\label{eqn:formula2}
\frac{d^2R}{dt^2}(s)=\mathrm{tr}\Big(A(s)^{-1}\frac{ d^2A}{dt^2}(s)\Big)-\mathrm{tr}\Bigg(\Big(A(s)^{-1}\frac{ dA}{dt}(s)\Big)^2\Bigg).
\end{equation}

We should also need the following version of the Whitney extension theorem (see \cite{Whitney} or \cite{Stein}).

\begin{theorem}
Let $A\subset\mathbb{R}^n$ be an arbitrary set, let $f\colon A\to\mathbb{R}$ and $V\colon A\to\mathbb{R}^n$. Suppose that there exists $M\in\mathbb{R}$ such that for all $x,y\in A$
\begin{align*}
&\abs{f(x)}\leq M, \norm{V(x)}\leq M,\\
& \norm{V(x)-V(y)}\leq M\norm{x-y},\\
&\abs{f(y)-f(x)-\langle V(x),y-x\rangle}\leq M\norm{x-y}^2.
\end{align*}
Then there exists a differentiable function $\tilde{f}\colon\mathbb{R}^n\to\mathbb{R}$ with locally Lipschitz derivative such that 
\begin{equation*}
\tilde{f}(x)=f(x), Df(x)(y)=\langle V(x),y\rangle\text{ for all }x\in A\text{ and all }y\in\mathbb{R}^n.
\end{equation*}
\end{theorem}

\begin{proof}[Proof of Theorem \ref{thm:discd}]
By Lemma \ref{lem:measurable} and Remark \ref{rem:cluster} it is enough to prove Theorem \ref{thm:discd} assuming that there is a single cluster of leaves $T_{sij}$.
Thus, let us fix a cluster $T_{sij}$.

Note that, by Corollary \ref{col:strength}, on $T_{sij}^{\lambda}$, $Du$ is Lipschitz; see Lemma \ref{lem:efge} for the definition of $T_{sij}^{\lambda}$. Moreover, by the second assertion of Lemma \ref{lem:important}, for any $x,y\in T_{sij}^{\lambda}$ there is
\begin{equation*}
\norm{u(y)-u(x)-Du(x)(y-x)}\leq \frac1{\lambda}\norm{x-y}^2.
\end{equation*}
By the Whitney extension theorem there exists a differentiable map $\tilde{u}$ with locally Lipschitz derivative on $\mathbb{R}^n$ that coincides with $u$ on $T_{sij}^{\lambda}$ and such that $D\tilde{u}=Du$ on $T_{sij}^{\lambda}$. By \cite[Lemma 3.2.4]{Klartag}, the second derivative of $\tilde{u}$ exists almost everywhere and is symmetric, in the sense that the second derivative of any of its components is symmetric. We will abuse the notation and assume that $u$ has Lipschitz derivative, is defined on $\mathbb{R}^n$, and its second derivative is symmetric $\lambda$-almost everywhere. 

Since $F\colon G(T_{sij}^0)\to\mathbb{R}^n$ has locally Lipschitz inverse, it follows that for $\lambda$-almost every $(a,b)\in G(T_{sij}^0)$ there exists $D^2u(F(a,b))$ and is symmetric.

By Fubini's theorem we infer that there exists a Borel cover  $(\Lambda_l)_{l=1}^{\infty}$ of $\Lambda$ such that for each $l\in\mathbb{N}$ there exists $b_l$ such that $(a,b_l)\in G(T_{sij}^0)$ for all $a\in\Lambda_l$. Moreover for $\lambda$-almost every $a\in \Lambda_l$, there exists $D^2u(F(a,b_l))$ and it is symmetric. Note that for $a\in\Lambda_l$
\begin{equation*}
u(F(a,b_l))=u(v(a))+b_l=s+b_l.
\end{equation*}
Hence, on the level set of $u$ corresponding to $s+b_l$, there exists $D^2u$ and it is symmetric.
Therefore, without loss of generality, passing to a refinement of initial cover and modifying the clusters $T_{sij}$, we assume that $D^2u(v(a))$ exists for $\lambda$-almost every $a\in\Lambda$ and it is symmetric.

Let $\mu$ have density $e^{-\rho}$ with respect to the Lebesgue measure on $\mathbb{R}^n$. For a leaf $\mathcal{S}$ such that $\mathrm{int}\mathcal{S}\subset T_{sij}^0$ and any Borel set $A\subset\mathbb{R}^n$ set
\begin{equation}\label{eqn:mudef}
\mu_{\mathcal{S}}(A)=\int_{\mathrm{int}\mathcal{S}}\mathbf{1}_A e^{-\rho}J_nF\circ G d\mathcal{H}_m.
\end{equation}
By Lemma \ref{lem:density} it follows now that for any Borel set $A\subset T_{sij}$ there is
\begin{equation*}
\mu(A)=\int_{\Lambda}\mu_{\mathcal{S}(v(a))}(A)d\lambda(a),
\end{equation*}
where $v\colon \mathbb{R}^{n-m}\to\mathbb{R}^n$ is the map from Lemma \ref{lem:cover}. Neglecting a set of Lebesgue measure zero, we may assume that $v$ is differentiable on the set $\Lambda$.
Let $\nu$ denote the push-forward of $\lambda$ via the map 
\begin{equation}\label{eqn:map}
\Lambda\ni a\mapsto \mathcal{S}(v(a))\in T^m.
\end{equation}
By Lemma \ref{lem:density} and the definition of $\nu$ the condition \ref{i:measur}) is satisfied.
Note that the map (\ref{eqn:map}) is Borel measurable, by the proof of Theorem \ref{thm:dis}. Hence $\nu$ is a Borel measure. For any Borel set $A\subset T_{sij}$ there is
\begin{equation*}
\mu(A)=\int_{T^m}\mu_{\mathcal{S}}(A)d\nu(\mathcal{S}).
\end{equation*}
Hence the condition \ref{i:disinteg}) of Theorem \ref{thm:discd} is satisfied. Condition \ref{i:concent}) holds true by the definition (\ref{eqn:mudef}). We shall prove that \ref{i:curv}) holds true as well.

Note that the density of a measure $\mu_{\mathcal{S}}$ for an $m$-dimensional leaf $\mathcal{S}$ is equal to
\begin{equation*}
\frac{d\mu_{\mathcal{S}}}{d\mathcal{H}_m}= J_nF\circ Ge^{-\rho}\mathbf{1}_{\mathcal{S}}.
\end{equation*}
Recall (see Lemma \ref{lem:efge}) that $F,G$ are given by the formulae
\begin{equation}\label{eqn:formu}
F(a,b)=v(a)+Du(v(a))^*(b)\text{ and }G(x)=(w(z),u(x)-u(z)),
\end{equation}
where $w\colon \mathbb{R}^n\to\mathbb{R}^{n-m}$ and $v\colon \mathbb{R}^{n-m}\to\mathbb{R}^n$ are maps from Lemma \ref{lem:cover}. Let us recall that $v(a)\in S_s^i$ for all $a\in\Lambda$. 
It follows by the definition of $S_s$ that $u(v(a))=s$ for all $a\in\Lambda$. Recall that, by Lemma \ref{lem:diff}, $u$ is differentiable in $T_{sij}^0$. Thus, as we assumed that $v$ is differentiable in $\Lambda$, for every $a\in\Lambda$ 
\begin{equation}\label{eqn:compo}
Du(v(a))Dv(a)=0.
\end{equation}
For $(a,b)\in G(T_{sij}^0)$ the derivative of $F$ at $(a,b)$ is equal to
\begin{equation*}
DF(a,b)=[Dv(a)+D^2u(v(a))^*(Dv(a)(\cdot))(b),Du(v(a))^*].
\end{equation*}
Note that for any vectors $z\in\mathbb{R}^{n-m}$ and $w\in\mathbb{R}^m$ the derivatives $Dv(a)z$ and $Du(v(a))^*w$ are orthogonal.
Indeed, by (\ref{eqn:compo}),
\begin{equation*}
\Big\langle Du(v(a))^*(w),Dv(a)(z)\Big\rangle=\Big\langle w, Du(v(a))Dv(a)(z)\Big\rangle=0.
\end{equation*}
Let $P$ denote the orthogonal projection onto the tangent space $V$ of the leaf $\mathcal{S}$ containing $v(a)$. Then by Lemma \ref{lem:diff}  $Du(v(a))=TP$.
Let $P^{\perp}$ denote the orthogonal projection onto the orthogonal complement of $V$. Then
\begin{equation*}
DF(a,b)=[Dv(a)+D^2u(v(a))^*(P^{\perp}Dv(a)(\cdot))(b),Du(v(a))^*].
\end{equation*}
Therefore, by the formula for block matrices, and as $Du(v(a))^*$ is isometric, we have
\begin{equation*}
\abs{\det(DF(a,b))}=\Big\lvert\det\Big(Dv(a)+P^{\perp}D^2u(v(a))^*(P^{\perp}Dv(a)(\cdot))(b)\Big)\Big\rvert,
\end{equation*}
which is equal to 
\begin{equation}\label{eqn:det}
\abs{\det \big(P^{\perp}Dv(a)\big)}\Big\lvert\det\Big(\mathrm{Id}+P^{\perp}D^2u(v(a))^*(P^{\perp}(\cdot))(b)\Big)\Big\rvert.
\end{equation}
Note that 
\begin{equation}\label{eqn:ha}
H(b)=\Big(\mathrm{Id} +P^{\perp}D^2u(v(a))^*(P^{\perp}(\cdot))(b)\Big)
\end{equation}
is a linear operator on the image of $P^{\perp}$, which is of dimension $n-m$. Moreover it is symmetric and invertible for any $b$ such that $(a,b)\in G(\mathrm{int}T_{pij})$, as $F$ is a bijection. Consider for some $b'\in\mathbb{R}^m$
\begin{equation*}
P^{\perp}D^2u(v(a))^*(P^{\perp}(\cdot))(b').
\end{equation*}
Let $A$ be such that
\begin{equation}\label{eqn:aa}
P^{\perp}D^2u(v(a))^*(P^{\perp}(\cdot))(b')=A\Big(\mathrm{Id} +P^{\perp}D^2u(v(a))^*(P^{\perp}(\cdot))(b)\Big).
\end{equation}
Then $A$ is conjugate to a symmetric operator of rank at most $n-m$, as
\begin{equation*}
H(b)^{-\frac12}AH(b)^{\frac12}=H(b)^{-\frac12}P^{\perp}D^2u(v(a))^*(P^{\perp}(\cdot))(b')H(b)^{-\frac12}.
\end{equation*}
In consequence, by the Cauchy-Schwarz inequality,
\begin{equation}\label{eqn:bac}
(\mathrm{tr}A)^2\leq (n-m)\mathrm{tr}(A)^2.
\end{equation}
Let $x=F(a,b)$. Let $q$ belong to the tangent space of $\mathcal{S}$. It is necessarily of the form $q=Du(v(a))^*(b')$ for some $b'\in\mathbb{R}^m$.  Then by (\ref{eqn:formu}), (\ref{eqn:det}) and (\ref{eqn:ha})
\begin{align*}
D\log\abs{\det DF\circ G}(x)(q)&=\frac{d}{dt}\log \abs{\det\big( DF(G(F(a,b)+tDu(v(a))^*(b')\big)}=\\
&=\frac{d}{dt}\log\abs{\det (DF(a,b+tb'))}=\frac{d}{dt}\log\abs{\det H(b+tb')}.
\end{align*}
Therefore, by (\ref{eqn:formula}), (\ref{eqn:formula2}) and by (\ref{eqn:aa})
\begin{equation*}
D\log\abs{\det DF\circ G}(x)(q)=\mathrm{tr}\Big(H(b)^{-1}P^{\perp}D^2u(v(a))^*(P^{\perp}(\cdot))(b')\Big)=\mathrm{tr}A
\end{equation*}
and
\begin{equation*}
D^2\log\abs{\det DF\circ G}(x)(q,q)=-\mathrm{tr}\Big(H(b)^{-1}P^{\perp}D^2u(v(a))^*(P^{\perp}(\cdot))(b')\Big)^2=-\mathrm{tr}(A^2).
\end{equation*}
By (\ref{eqn:bac}) and by Lemma \ref{lem:tri}, if $N\notin [m, n]$, then
\begin{align*}
-D^2 \log\abs{\det DF\circ G}(x)(q,q)&=\mathrm{tr}(A^2)\geq \\
&\geq\frac{1}{n-m}(\mathrm{tr}A)^2\geq \frac1{N-m}(D\rho(x) (q)-\mathrm{tr} A)^2-\frac{(D\rho(x)(q))^2}{N-n}.
\end{align*}
Note that by the assumption on $\mu$, c.f. Definition \ref{defin:curv}, for all $p\in\mathbb{R}^n$
\begin{equation*}
D^2{\rho}(x)(p,p)-\frac{D\rho(x)(p)^2}{N-n}\geq \kappa \norm{p}^2.
\end{equation*}
Thus for all $q$ in the tangent space of $\mathcal{S}$ there is
\begin{equation*}
D^2\rho(x)(q,q)-D^2 \log\abs{\det DF\circ G}(x)(q,q)-\frac{\big(D\rho(x)(q)-D(\log\abs{\det DF\circ G})(x)(q)\big)^2}{N-m}\geq \kappa \norm{q}^2.
\end{equation*}
We infer that $(\mathrm{int}\mathcal{S},\norm{\cdot},\mu_{\mathcal{S}})$ satisfies the curvature-dimension condition $CD(\kappa,N)$, provided that $N\notin[m,n]$.

If $N=n$, then $\rho$ is required to be a constant function, and thus in this case the inequality is also satisfied. If $N=\infty$, then the desired estimates follow readily.
\end{proof}
For the historical remarks on similar estimates we refer to \cite{Klartag}. 
\section*{Acknowledgements}
The author wishes to thank Bo'az Klartag for proposing to work on this problem and for useful discussions.

  \bibliographystyle{plain} 
  \bibliography{refs}

\end{document}